\documentclass[12pt]{amsart}
\usepackage{amsmath,amsfonts,amsthm,amssymb,amscd,stmaryrd,url}

\newcommand{\st}{s.t.}

\DeclareFontFamily{OT1}{rsfs}{}
\DeclareFontShape{OT1}{rsfs}{n}{it}{<-> rsfs10}{}
\DeclareMathAlphabet{\mathscr}{OT1}{rsfs}{n}{it}
\DeclareMathOperator{\sgn}{sgn}

\DeclareMathOperator{\num}{num}
\DeclareMathOperator{\den}{den}

\DeclareMathOperator{\rad}{rad}
\DeclareMathOperator{\meas}{meas}
\DeclareMathOperator{\mo}{\,mod}
\DeclareMathOperator{\sq}{sq}

\DeclareMathOperator{\Area}{Area}

\DeclareMathOperator{\Disc}{Disc}

\DeclareMathOperator{\irr}{irr}
\DeclareMathOperator{\rnk}{rank}

\DeclareMathOperator{\Frob}{Frob}
\DeclareMathOperator{\Gal}{Gal}
\DeclareMathOperator{\Cl}{Cl}

\newtheorem{prop}{Proposition}[section]
\newtheorem{thm}[prop]{Theorem}

\newtheorem{cor}[prop]{Corollary}
\newtheorem{lem}[prop]{Lemma}
\newtheorem{defn}{Definition}
\newtheorem{Exa}{Example}
\newtheorem*{defn*}{Definition}
\newenvironment{Rem}{{\bf Remark.}}{}
\numberwithin{equation}{section}
\author{H. A. Helfgott}
\address{H. A. Helfgott, D\'epartement de Math\'ematiques et Statistique, Universit\'e de Montr\'eal, Montr\'eal,
QC H3C 3J7, Canada}
\email{helfgott@dms.umontreal.ca}
\title{On the square-free sieve}
\begin{document}
\maketitle

\section{Introduction}
A \emph{square-free
sieve} is a result that gives an upper bound for how often a square-free 
polynomial may adopt values that are not square-free.
More generally, we may wish to 
control the behavior of a function depending on the largest square factor of $P(x_1,\dotsc,x_n)$, where $P$ is a square-free polynomial.

We may aim at obtaining an asymptotic expression
\begin{equation}\label{eq:bartok}
\text{main term } + O(\text{error term}),
\end{equation}
where the main term will depend on the application; in general, the
error term will depend only on the polynomial $P$ in question,
not on the particular quantity being estimated. We can split the
error term further into one term that can be bounded for
every given $P$, and a second term, say, $\delta(P)$, which may be
rather hard to estimate, and which is unknown for polynomials $P$
of high enough degree. Given this framework, the strongest
results in the literature may be summarized as follows:

\begin{center} \begin{tabular}{l|ll}
$\deg_{\irr}(P)$ &$\delta(P(x))$ &$\delta(P(x,y))$\\ \hline
$1$ & $\sqrt{N}$ & $1$\\
$2$ & $N^{2/3}$ & $N$\\
$3$ & $N / (\log N)^{1/2}$ & $N^2 / \log N$ \\
$4$ & & $N^2/\log N$ \\
$5$ & & $N^2/\log N$ \\
$6$ & & $N^2/(\log N)^{1/2}$
\end{tabular} \end{center}

Here $\deg_{\irr}(P)$ denotes the degree of the largest irreducible
factor of $P$. The second column gives $\delta(P)$ for polynomials
$P\in \mathbb{Z}\lbrack x\rbrack$ of given $\deg_{\irr}(P)$,
whereas the third column refers to homogeneous polynomials 
$P\in \mathbb{Z}\lbrack x,y\rbrack$. The trivial estimates would
be $\delta(P(x))\leq N$ and $\delta(P(x,y))\leq N^2$. See section
\ref{sec:bestpast}
for attributions.

Our task can be divided into two halves. The first one, undertaken in 
section \ref{sec:cuacriba}, consists in estimating all terms but
$\delta(N)$. We do as much in full generality for any $P$ over
any number field. The second half regards bounding
$\delta(N)$. We improve on all estimates known for $\deg P\geq 3$:

\begin{center}\begin{tabular}{l|ll}
$\deg_{\irr}(P)$ &$\delta(P(x))$ &$\delta(P(x,y))$\\ \hline 
$3$ & $N/(\log N)^{0.6829\dotsc}$ & $N^{3/2}/\log N$\\
$4$ & & $N^{4/3} (\log N)^A$ \\
$5$ & & $N^{(5 + \sqrt{113})/8+\epsilon}$\\
$6$ & & $N/(\log N)^{0.7043\dotsc}$
\end{tabular}\end{center}

The bound for $\deg_{\irr}(P) = 6$ depends on the Galois group of $P(x,1)$;
the table gives the bound for the generic case $\Gal = S_6$. Some important
applications, e.g., the computation of the distribution of root numbers
in generic families of elliptic curves, yield an expression such as 
(\ref{eq:bartok}) with the main term equal to zero. Our improvements on
the error term then become (substantial)
improvements on the upper and lower bounds.

Most of our improvements hinge on a change from a local to a global
perspective. Such previous work in the field as was purely
sieve-based can be seen as an series of purely local estimates
on the density of points on curves of non-zero genus. Our techniques
involve a mixture of sieves, height functions and
sphere packings.

The present paper is motivated in part by applications
that demand us to construe square-free sieves in a broader sense than may
have been customary before. We sketch the most general
setting in which a square-free sieve is meaningful and applicable,
and provide statements of intermediate generality and concreteness.
In particular, we show how to use a square-free sieve to 
give explicit estimates on the average of certain infinite products
arising from $L$-functions.

\section{Notation}
Let $n$ be a non-zero integer.
We write $\tau (n)$ for the number of positive 
divisors of $n$, $\omega (n)$ for the
number of the prime divisors of $n$, and $\rad(n)$ for the product of the prime divisors
of $n$. Given a prime $p$, we write $v_p(n)$ for the non-negative integer $j$
such that $p^j|n$ and $p^{j+1}\nmid n$. 
For any $k\geq 2$, we write $\tau_k(n)$ for the number of 
$k$-tuples $(n_1,n_2,\dotsc, n_k)\in (\mathbb{Z}^+)^k$ such that
$n_1\cdot n_2\cdot \dotsb n_k = |n|$. Thus $\tau_2(n) = \tau(n)$.
We adopt the convention that $\tau_1(n) = 1$. 
We let
\[\sq(n) = \prod_{p^2|n} p^{v_p(n)-1} .\]
We call a rational integer $n$ {\em square-full} if $p^2|n$ for every prime $p$
dividing $n$. Given any non-zero rational integer $D$, we say that $n$ is
{\em ($D\!$)-square-full} if $p^2|n$ for every prime $p$ that divides $n$ but not $D$.

We denote by $\mathscr{O}_K$ the ring of integers of a global or local
field $K$. We let $I_K$ be the semigroup of non-zero ideals of $\mathscr{O}_K$.
Given a non-zero ideal $\mathfrak{a}\in I_K$, 
we write $\tau_K(\mathfrak{a})$ for the number of ideals dividing
 $\mathfrak{a}$,
$\omega_K(\mathfrak{a})$ for the number of prime ideals dividing 
$\mathfrak{a}$,
and $\rad_K(\mathfrak{a})$ for the product of the prime ideals dividing
$\mathfrak{a}$.
Given a positive integer $k$, we write $\tau_{K,k}(\mathfrak{a})$ for the number of 
$k$-tuples $(\mathfrak{a}_1,\mathfrak{a}_2,\dotsc ,\mathfrak{a}_k)$ of
ideals of $\mathscr{O}_K$ such that
$\mathfrak{a} = \mathfrak{a}_1 \mathfrak{a}_2 \dotsb \mathfrak{a}_k$. 
Thus $\tau_2(\mathfrak{a}) = \tau(\mathfrak{a})$.
We define 
$\rho(\mathfrak{a})$ to be the positive integer generating $\mathfrak{a}\cap \mathbb{Z}$.

When we say that a polynomial $f\in \mathscr{O}_K\lbrack x\rbrack$ or
$f\in K\lbrack x\rbrack$ is {\em square-free}, we always mean that
$f$ is square-free as an element of $K\lbrack x\rbrack$. Thus, for
example, we say that $f\in \mathbb{Z}\lbrack x\rbrack$ is {\em square-free}
if there is no polynomial $g\in \mathbb{Z}\lbrack x\rbrack$
such that $\deg g \geq 1$ and $g^2|f$. 

Given an elliptic curve $E$ over $\mathbb{Q}$, we write $E(\mathbb{Q})$
for the set of rational (that is, $\mathbb{Q}$-valued) points of $E$.
We denote by $\rnk(E)$ the algebraic rank of $E(\mathbb{Q})$.

We write $\# S$ for the cardinality of a finite set $S$. 
Given two sets $S_1\subset
S_2$, we denote $\{x \in S_2 : x\notin S_1\}$ by $S_2 - S_1$.

By a {\em lattice} we will always mean an additive subgroup of 
$\mathbb{Z}^2$ of 
finite index. For $S\subset \lbrack -N, N\rbrack$ a convex set and $L$
a lattice,
\[\#(S \cap L) = \frac{\Area(S)}{\lbrack \mathbb{Z}^2 : L \rbrack} + O(N),\]
where the implied constant is absolute. 

A {\em sector} is a connected component of a set of the form
$\mathbb{R}^2 - (T_1 \cup \dotsb \cup T_n)$, where $n$ is a non-negative
integer and $T_i$ is a line going
through the origin. Every sector $S$ is convex.
\section{Sieving}\label{sec:cuacriba}
\subsection{Averaging infinite products}
Take a function $u:\mathbb{Z}\to \mathbb{C}$ defined as a product
\[u(n) = \prod_p u_p(n) \]
of functions $u_p:\mathbb{Z}_p \to \mathbb{C}$, one for each prime.
One may wish to write the average of $u(n)$ as a product of $p$-adic
integrals:
\begin{equation}\label{eq:methy}
\lim_{N\to \infty} \frac{1}{N} \sum_{n=1}^N u(n) = \prod_p
 \int_{\mathbb{Z}_p} u_p(x) dx .\end{equation}
 Unfortunately, such an expression is not valid in general. Even if all 
the integrals are defined and of modulus at most $1$,
the infinite product may not converge, and the average on the left
of (\ref{eq:methy}) may not be defined; even if the product does 
converge and the average is defined, the two may not be equal. (Take,
for example, $u_p(x) = 1$ for $x\in \mathbb{Z}$, $u_p(x) = 0$ for
$x\in \mathbb{Z}_p - \mathbb{Z}$.)

We will establish that (\ref{eq:methy}) is in fact true when three conditions
hold. The first one will be a local condition ensuring  that
\[\lim_{N\to \infty} \frac{1}{N} \sum_{n=1}^N u_p(n) = 
\int_{\mathbb{Z}_p} u_p(x) dx,\]
among other things. The second condition states that $u_p(x)=1$ 
when $p^2\nmid P(x)$, where $P(x)$
is a fixed polynomial with integer coefficients. 
Thus we know that we are being asked to evaluation a convergent product
such as $\prod_p (1 - a_p/p^2)$, rather than, say, the product
$\prod_p (1 - 1/p)$. The third and last condition is that there be
a non-trivial bound on the term $\delta(P(x))$ mentioned in the introduction.

We can make classical square-free sieves fit into the present framework
by setting $u_p(x)=0$ when $p^2|F(x)$ (see section \ref{sec:cuanumer}).
In other instances, the full generality of our treatment becomes 
necessary: see \cite{He} for the example
\begin{equation}\label{eq:agort}
u(n) = \prod_p u_p(n) = \prod_p \frac{W_p(\mathscr{E}(n))}{W_{0,p}(
\mathscr{E}(n))},\end{equation}
where $\mathscr{E}$ is a family of elliptic curves, $W_p(E)$ is
the local root number of an elliptic curve $E$, and $W_{0,p}(E)$ is
a first-order approximation to $W_p(E)$. It seems reasonable to expect
the present method to be applicable to the estimation of other such
ratios arising from Euler products.
\begin{defn}\label{def:1900}
A function $f:\mathbb{Z}_p \to \mathbb{C}$ is 
\emph{openly measurable} if, for every $\epsilon>0$, there is a 
partition $\mathbb{Z}_p = U_0 \cup U_1 \cup \dotsc \cup 
U_{k}$ and a tuple of complex numbers 
$(y_1,\dotsc , y_{k})$ such that
\begin{enumerate}
\item $|f(x) - y_j| < \epsilon$ for $x\in U_j$, $1\leq j\leq k$,
\item $\meas(U_0) < \epsilon$,
\item all $U_j$ are open.
\end{enumerate}
\end{defn}
Every function $f:\mathbb{Z}_p \to \mathbb{C}$ continuous outside
a set of measure zero is openly measurable.
\begin{prop}\label{prop:unfair}
For every prime $p$, let
$u_p: \mathbb{Z}_p \to \mathbb{C}$ be an openly
measurable function with $|u_{\mathfrak{p}}(x)|\leq 1$ for every
$x\in \mathbb{Z}_p$. Assume that $u_p(x)=1$
unless $p^2 |P(x)$, where $P\in \mathbb{Z}\lbrack x \rbrack$
is a polynomial satisfying
\[\# \{ 1\leq x \leq N : \exists p>N^{1/2},
p^2 | P(x,y)\} = o(N) .\]
Let  $u(n) = \prod_{p} u_{p}(n)$. Then (\ref{eq:methy}) holds.
\end{prop}
\begin{proof}
Since we are not yet being asked to produce estimates for the speed of 
convergence, we may adopt a simple procedure specializing to the ones
in \cite{Ho}, Ch. IV, and \cite{Gr}. Let $\epsilon_0$ be given.
By Appendix \ref{sec:appa}, Lemma \ref{lem:sols},
\[\left| 1 - \int_{\mathbb{Z}_p} u_p(x) d x \right| = O(p^{-2}) ,\]
where the implied constant depends only on $P$. Hence
\[\left|\prod_p \int_{\mathbb{Z}_p} u_p(x) dx -
\prod_{p\leq \epsilon_0^{-1}} \int_{\mathbb{Z}_p} u_p(x) dx \right| = 
O(\epsilon_0) .\]
Let \[u'(x) = \prod_{p\leq \epsilon_0^{-1}} u_p(x).\]
Since $u'(x)=u(x)$ unless $p^2|P(x)$ for some $p > \epsilon_0^{-1}$,
\[\begin{aligned}
\left|\sum_{n=1}^N u(n) - \sum_{n=1}^N u'(n) \right| &\leq
2 \cdot \# \{1\leq x\leq N : \exists p > \epsilon_0^{-1} \,\st\, p^2|P(x)\} \\
&\leq 2 \sum_{\epsilon_0^{-1} < p \leq N^{1/2}} 
 \mathop{\sum_{1\leq x\leq N}}_{p^2|P(x)} 1\\
&+ 2 \cdot \# \{1\leq x\leq N : \exists p > N^{1/2}, p^2|P(x)\}.\end{aligned}\]
By Appendix \ref{sec:appa}, Lemma \ref{lem:dontaskmemyadvice}, 
\[\sum_{\epsilon_0^{-1} < p \leq N^{1/2}} 
\mathop{\sum_{1\leq x \leq N}}_{p^2 | P(x)} 1 =
O\left(\sum_{\epsilon_0^{-1} < p \leq N^{1/2}} N/p^2\right) 
= O(\epsilon_0 N) ,\]
where the implied constant depends only on $P$. By the assumption in
the statement,
\[\# \{1\leq x \leq N : \exists p > N^{1/2}, p^2 | P(x)\}\leq \epsilon_0 N\]
for $N$ greater than some constant $C_{\epsilon_0}$.
It now remains to compare
\[\frac{1}{N} \sum_{n=1}^N u'(n) \text{ \; and \; }
  \prod_{p<\epsilon_0^{-1}} \int_{\mathbb{Z}_p} u_p(x) d x .\]
For every $p$, we are given a tuple $(y_1,\dotsc,y_{k_p})$ and a partition
$\mathbb{Z}_p = U_{p,0} \cup U_{p,1} \cup \dotsb \cup U_{p,k_p}$ 
satisfying the conditions in Definition \ref{def:1900} with 
$\epsilon = \epsilon_0^2$. We can assume that $U_{p,1},U_{p,2},\dotsc,
U_{p,k_p}$ are connected and $|y_1|,\dotsc,|y_{k_p}|\leq 1$. We have
\[\begin{aligned}
\sum_{n=1}^N u'(n) &= \sum_{\vec{j} \in \prod_{p<\epsilon_0^{-1}}
\lbrack 0,k_p \rbrack} \mathop{\sum_{1\leq n\leq N}}_{n\in \bigcap_p
 U_{p,j_p} \cap \mathbb{Z}} u'(n) \\
&= \sum_{\vec{j} \in \prod_{p<\epsilon_0^{-1}} \lbrack 0,k_p\rbrack}
 \mathop{\sum_{1\leq n\leq N}}_{n\in \bigcap_p U_{p,j_p}\cap \mathbb{Z}}
\left(\prod_{p<\epsilon_0^{-1}} y_{p,j_p} + 
O\left((1 + \epsilon_0^2)^{\epsilon_0^{-1}} - 1\right)\right) \\
&= \sum_{\vec{j} \in \prod_{p<\epsilon_0^{-1}} \lbrack 0,k_p\rbrack}
 \mathop{\sum_{1\leq n\leq N}}_{n\in \bigcap_p U_{p,j_p}\cap \mathbb{Z}}
 \prod_{p<\epsilon_0^{-1}} y_{p,j_p} + O(\epsilon_0 N) \\
&= \sum_{\vec{j} \in \prod_{p<\epsilon_0^{-1}} \lbrack 0,k_p\rbrack}
\prod_{p<\epsilon_0^{-1}} y_{p,j_p} \cdot \#\{1\leq n\leq N : n \in
 \bigcap_p U_{p,j_p} \cap \mathbb{Z}\} \\ &+ O(\epsilon_0 N),\end{aligned}\]
whereas
\[\begin{aligned}
\prod_{p<\epsilon_0^{-1}} \int_{\mathbb{Z}_p} u_p(x) d x \,&=
\sum_{\vec{j} \in \prod_{p<\epsilon_0^{-1}} \lbrack 0, k_p \rbrack}\, 
 \prod_{p<\epsilon_0^{-1}} \int_{U_{p,j_p}} u_p(x) d x \\
&= \sum_{\vec{j} \in \prod_{p<\epsilon_0^{-1}} \lbrack 0, k_p \rbrack}\, 
 \prod_{p<\epsilon_0^{-1}} \meas(U_{p,j_p}) (y_{p,j_p} + O(\epsilon)) .
\end{aligned}\]
The absolute value of
\[
\sum_{\vec{j} \in \prod_{p<\epsilon_0^{-1}} \lbrack 0, k_p \rbrack}\, 
 \prod_{p<\epsilon_0^{-1}} \meas(U_{p,j_p}) (y_p + O(\epsilon)) -
\sum_{\vec{j} \in \prod_{p<\epsilon_0^{-1}} \lbrack 0, k_p \rbrack}\, 
 \prod_{p<\epsilon_0^{-1}} \meas(U_{p,j_p}) y_p\]
is at most a constant times 
\[\begin{aligned}
\sum_{\vec{j} \in \prod_{p<\epsilon_0^{-1}} \lbrack 0, k_p \rbrack}
 &\left(\prod_{p<\epsilon_0^{-1}} \meas(U_{p,j_p}) ( 1 + \epsilon) 
 - \prod_{p<\epsilon_0^{-1}} \meas(U_{p,j_p})\right) 
\\ &\ll \prod_{p<\epsilon_0^{-1}} (1 + \epsilon) - 1 \ll \epsilon_0 .
\end{aligned}\]
It is left to bound the difference between
\[\frac{1}{N} \sum_{\vec{j} \in \prod_{p<\epsilon_0^{-1}} 
\lbrack 0, k_p\rbrack } y_{p,j_p}\cdot
\# \{1\leq n\leq N : n\in \bigcap_p U_{p,j_p} \cap \mathbb{Z}\}\]
and
\[\sum_{\vec{j} \in \prod_{p<\epsilon_0^{-1}} 
\lbrack 0, k_p\rbrack } \, \prod_{p<\epsilon_0^{-1}} \meas(U_{p,j_p}) y_{p,j_p}
.\]
It is enough to estimate
\begin{equation}\label{eq:spikop}
\sum_{\vec{j} \in \prod_{p<\epsilon_0^{-1}} 
\lbrack 0, k_p\rbrack } \left(\prod_{p<\epsilon_0^{-1}} \meas(U_{p,j_p}) 
- \frac{1}{N} 
\# \{1\leq n\leq N : n\in \bigcap_p U_{p,j_p} \cap \mathbb{Z}\}\right) .
\end{equation}
When all $j_p$ are positive, $\bigcap_p U_{p,j_p} \cap \mathbb{Z}$ is
an arithmetic progression. Hence
\[
N \prod_j \meas(U_{p,j_p}) -
\#\{1\leq n\leq N : n \in \bigcap_p U_{p,j_p} \cap \mathbb{Z}\} \leq 1.\]
When some $j_p$ are zero, we can use inclusion-exclusion to obtain
\[\left|\#\{1\leq n\leq N : n \in \bigcap_p U_{p,j_p} \cap \mathbb{Z}\} 
- N \prod_j \meas(U_{p,j_p})\right| \leq
\mathop{\prod_p}_{j_p=0} k_p .\]
Hence the absolute value of (\ref{eq:spikop}) is at most
\[\frac{1}{N} \prod_{p<\epsilon_0^{-1}} (2 k_p) .\]
If $N>\epsilon_0^{-1} \prod_{p<\epsilon_0^{-1}} (2 k_p)$, then clearly
$\frac{1}{N} \prod_{p<\epsilon_0^{-1}} (2 k_p) = O(\epsilon_0)$.

We can conclude that, for $N\geq \max(C_{\epsilon_0}, \epsilon_0^{-1}
\prod_{p<\epsilon_0^{-1}} (2 k_p))$, 
\[\left|\lim_{N\to \infty} \sum_{n=1}^N u(n) - 
\prod_p \int_{\mathbb{Z}_p} u_p(x) dx\right| = O(\epsilon_0),\]
where the implied constant depends only on $P$. 
\end{proof}
The question now is to make Prop. \ref{prop:unfair} explicit, or, rather,
how to do so well. We desire strong bounds on the speed of convergence.
\subsection{Riddles}
We will now see a near-optimal
way to sieve out square factors below a certain
size. Since the procedure is highly formal, we will state it in fairly general
terms. It will be effortless to derive statements on number fields and 
higher power divisors.

We write $P(\mathscr{P})$ for the set of all subsets of a given set
$\mathscr{P}$.
\begin{defn}
A {\em soil} is a tuple $(\mathscr{P},\mathscr{A},r,f)$, where $\mathscr{P}$
is a set, $\mathscr{A}$ is a finite set, $r$ is a function from $\mathscr{A}$
to $P(\mathscr{P})$, and $f$ is a function from 
$\mathscr{A} \times P(\mathscr{P})$ to $\mathbb{C}$.
\end{defn}
The purpose of a sieve is to estimate
\begin{equation}\label{eq:fireq}
\sum_{a\in \mathscr{A}} f(a,r(a))
\end{equation}
given data on
\begin{equation}\label{eq:seceq}
\mathop{\sum_{a\in \mathscr{A}}}_{r(a)\supset d_1} f(a,d_2)
\end{equation}
for $d_1, d_2 \in P(\mathscr{P})$. A traditional formulation would set
$f(a,d)=0$ for $d$ non-empty and $f(a,\emptyset)=1$ for every $a$. We shall work
with $f(a,d)$ bounded for the sake of simplicity.

We need a way to order $\mathscr{P}$. We will take as given a function
$h:P(\mathscr{P}) \to \mathbb{Z}^+$ such that
\begin{equation}
\renewcommand{\theequation}{h1}
\text{$h(d_1 \cup d_2) \leq h(d_1) h(d_2)$
for all $d_1, d_2 \in P(\mathscr{P})$ disjoint.}
\end{equation}
\addtocounter{equation}{-1}
and
\begin{equation}
\renewcommand{\theequation}{h2}
\text{$\{d\in P(\mathscr{P}) : h(d) \leq n\}$ is finite for
every $n\in \mathbb{Z}$.} \end{equation}
\addtocounter{equation}{-1}

We will also need an estimate for (\ref{eq:seceq}) in terms
of $h(d)$. In our applications, we will be able to assume
\begin{equation}
\renewcommand{\theequation}{A1}
\mathop{\sum_{a\in \mathscr{A}}}_{r(a)\supset d} 1 \leq C_0 
\frac{X C_1^{\# d}}{h(d)} + C_0 C_2^{\# d} 
\text{ \; for $d\subset \mathscr{P}$,}
\end{equation}
\addtocounter{equation}{-1}
\begin{equation}
\renewcommand{\theequation}{A2}
\mathop{\sum_{a\in \mathscr{A}}}_{r(a)\supset d_1} f(a,d_2) =
X \frac{g(d_1,d_2)}{h(d_1)} + r_{d_1,d_2}
\text{ \; for $d_2\subset d_1 \subset 
\mathscr{P}$, $h(d)\leq M_0$,}
\end{equation}
\addtocounter{equation}{-1}
where $X$, $C_0$, $C_1$ and $C_2$ are constants given by the soil, $g$ is some
function with desirable properties, and $r_{d_1,d_2}$ is small in average.
We shall mention explicitly when and whether we assume (h1),
(h2), (A1) and (A2) to hold; we will also state our precise 
conditions on $g$ and $r_d$ when we assume A2. We will henceforth write
$S_d$ and $A_{d_1,d_2}$ for the sums on the left of (A1) and (A2),
respectively.

The sieve we are about to propose is of use when $\{h(d)\}$ is fairly
sparse; hence the title.

We write $\mu(S)$ for $(-1)^{\# S}$.
\begin{prop}\label{prop:ridd}
Let $(\mathscr{P}, \mathscr{A}, r, f)$ be a soil with $f$ bounded. 
Let $h:P(\mathscr{P}) \to \mathbb{Z}^+$ satisfy (h1) and 
(h2). Then, for every positive integer $M$,
\[\left|\sum_{a\in \mathscr{A}} f(a,r(a)) -
 \mathop{\sum_{d\subset \mathscr{P}}}_{h(d)\leq M} 
\sum_{d'\subset d} \mu(d-d') A_{d,d'} \right|
\]
is at most
\[
\left(\mathop{\sum_{d\subset \mathscr{P}}}_{
M < h(d) \leq M^2} (3^{\# d} + 3) S_d +
\mathop{\sum_{p\in \mathscr{P}}}_{h(\{p\})>M^2} S_{\{p\}}\right) 
\cdot \max_{a,d} f(a,d) .\]
\end{prop}
\begin{proof}
For every $d\subset P$, let $\pi(d) = \{x\in d: h(\{x\}) \leq M\}$. By
M\"obius inversion,
\[\mathop{\sum_{d\subset r(a)}}_{x\in d \Rightarrow h(\{x\}) \leq M}
\sum_{d'\subset d}
\mu(d-d') f(a,d') = f(a,\pi(r(a))) \]
for every $a\in \mathscr{A}$. Hence
\[\begin{aligned}
\sum_{a\in \mathscr{A}} f(a,r(a)) &= \sum_{a\in \mathscr{A}}
 (f(a,r(a)) - f(a,\pi(r(a)))) + \sum_{a\in \mathscr{A}} \delta_a \\
&+ \mathop{\sum_{d\subset \mathscr{P}}}_{h(d) \leq M} \sum_{d'\subset d}
 \mu(d-d') A_{d,d'},\end{aligned}\]
where we write
\begin{equation}\label{eq:delt}
\delta_a = 
\mathop{\sum_{d\subset r(a)}}_{x\in d \Rightarrow h(\{x\}) \leq M}
\sum_{d'\subset d}
\mu(d-d') f(a,d') -
\mathop{\sum_{d\subset r(a)}}_{h(d) \leq M} \sum_{d'\subset d}
\mu(d-d') f(a,d') .\end{equation}
Since $a=\pi(a)$ unless some $x\in a$ satisfies $h(\{x\})>M$, we know that
\[\sum_{a\in \mathscr{A}}
 (f(a,r(a)) - f(a,\pi(r(a)))) \leq
2 \max_{a,d} f(a,d) \mathop{\sum_{x\in \mathscr{P}}}_{h(\{x\})>M}
S_{\{x\}} .\]
Now take $a\in \mathscr{A}$ such that $\delta_{a}\neq 0$. Then 
$h(\pi(r(a)))>M$. Let $d$ be a subset of $a$ with $h(d)\leq M$. We would
like to show that there is a subset $s$ of $r(a)$ such that $d\subset s$
and $M<h(s)\leq M^2$. Since $h(d)\leq M$, all elements $x\in d$ obey
$h(\{x\})\leq M$, and thus $d\subset \pi(r(a))$. Let $x_1,\dotsc,x_k$
be the elements of $\pi(r(a)) - d$. Let $s_0 = d$. For $1\leq i\leq k$,
let $s_i = d \cup \{x_1,\dotsc,x_i\}$. Then $h(s_0)\leq M$,
$h(s_k) = h(\pi(r(a))) > M$ and $h(s_{i+1}) \leq h(s_i) h(\{x_i\})
\leq h(s_i) \cdot M$ by (h1). Hence there is an $0\leq i<k$ such that
$M<h(s_i)\leq M^2$. Since $d\subset s_i$ and $s_i\subset \pi(r(a))$,
we can set $s = s_i$.

Now we bound the second sum in (\ref{eq:delt}) trivially:
\[\left|\mathop{\sum_{d\subset r(a)}}_{h(d)\leq M} \sum_{d'\subset d}
\mu(d-d') f(a,d') \right| \leq \max_{a,d} |f(a,d)| \cdot \mathop{\sum_{d\subset r(a)}}_{
h(d)\leq M} 2^{\# d}.\]
By the foregoing discussion,
\[\mathop{\sum_{d\subset r(a)}}_{h(d)\leq M} 2^{\# d} \leq
\mathop{\sum_{s\subset r(a)}}_{M< h(s) \leq M^2} \sum_{d\subset s} 2^{\# d}
= \mathop{\sum_{s\subset r(a)}}_{M<h(s) \leq M^2} 3^{\# s} .\]
(We are still assuming $\delta_a \ne 0$.) Since
\[\left|
\mathop{\sum_{d\subset r(a)}}_{x\in d \Rightarrow h(\{x\}) \leq M}
\sum_{d'\subset d}
\mu(d-d') f(a,d') \right| = |f(a,\pi(r(a)))| \leq \max_{a,d} f(a,d),\]
and since (again by $\delta_a\ne 0$) we have
\[\mathop{\sum_{s\subset r(a)}}_{M<h(s)\leq M^2} 1 \geq 1,\]
we can conclude that 
\[\sum_{a\in \mathscr{A}} \delta_a \leq \max_{a,d} |f(a,d)| \cdot
\sum_{a\in \mathscr{A}} \mathop{\sum_{s\subset r(a)}}_{M<h(s)\leq M^2}
 (3^{\# s} + 1) . \]
Clearly
\[\sum_{a\in \mathscr{A}} \mathop{\sum_{s\subset r(a)}}_{M<h(s)\leq M^2}
 (3^{\# s} + 1) \,\leq \mathop{\sum_{d\subset \mathscr{P}}}_{M < h(d) \leq M^2}
(3^{\# d} + 1) S_d .\]
The statement follows.
\end{proof}

\begin{cor}\label{cor:yugo}
Let $(\mathscr{P},\mathscr{A},r,f)$ be a soil with $\max_{a,d} f(a,d) \leq
C_3$. Assume (A1) and (A2) with $\max_{d_1,d_2} g(d_1,d_2)
\leq C_4$. Let h be multiplicative and satisfy (h2). Then, for every 
$M\leq M_0$,
\[\left|\sum_{a\in \mathscr{A}} f(a,r(a)) -
X\cdot 
\sum_{d\subset \mathscr{P}} \sum_{d'\subset d} \mu(d-d') \frac{g(d,d')}{
h(d)} 
\right|\]
is at most
\begin{equation}\label{eq:vivsav}\begin{aligned}
X &\mathop{\sum_{d\subset \mathscr{P}}}_{h(d)>M}
 \frac{1}{h(d)} (C_4 2^{\# d} + C_3 C_0 C_1^{\# d} (3^{\# d} + 3)) \\
&+ C_3\cdot \mathop{\sum_{d\subset \mathscr{P}}}_{M<h(d)\leq M^2} 
C_0 C_2^{\# d} (3^{\# d} + 3) \\
&+ \mathop{\sum_{d\subset \mathscr{P}}}_{h(d)\leq M}
 \sum_{d'\subset d} |r_{d,d'}| + C_3 \mathop{\sum_{p\in \mathscr{P}}}_{
h(\{p\}) > M^2} S_{\{p\}} .\end{aligned}\end{equation}
\end{cor}
\begin{proof}
Apply Proposition \ref{prop:ridd}. By (A1),
\[\begin{aligned}
\mathop{\sum_{d\subset \mathscr{P}}}_{M<h(d)\leq M^2} 
(3^{\# d} + 3) S_d
&\leq X \cdot \mathop{\sum_{d\subset \mathscr{P}}}_{M<h(d)\leq M^2} 
C_0 \frac{C_1^{\# d}}{h(d)} (
3^{\# d} + 3) \\ &+ \mathop{\sum_{d\subset \mathscr{P}}}_{M<h(d)\leq M^2}
C_0 C_2^{\# d} (3^{\# d} + 3) .\end{aligned}\]
By (A2),
\[\begin{aligned}
\mathop{\sum_{d\subset \mathscr{P}}}_{h(d)\leq M} \sum_{d'\subset d}
\mu(d-d') A_{d,d'} &= X \cdot \mathop{\sum_{d\subset \mathscr{P}}}_{h(d)\leq M}
\sum_{d'\subset d} \mu(d-d') \frac{g(d,d')}{h(d)} \\
&+ \mathop{\sum_{d\subset \mathscr{P}}}_{h(d)\leq M} \sum_{d'\subset d}
\mu(d-d') r_{d,d'} .\end{aligned}\]
Finally,
\[\left|\mathop{\sum_{d\subset \mathscr{P}}}_{h(d)>M} \sum_{d'\subset d}
\mu(d-d') \frac{g(d,d')}{h(d)} \right| \leq \mathop{\sum_{d\subset P}}_{
h(d)> M} \frac{C_4 2^{\# d}}{h(d)} .\]
\end{proof}
\begin{prop}\label{prop:owerr}
Let $K$ be a number field. Let $Q\in \mathscr{O}_K\lbrack x\rbrack$
be a polynomial. Let $m$ be a positive integer. Then the number of
positive integers $x\leq N$ for which $Q(x)$ is free of $m$th powers equals
\[\begin{aligned}
N \prod_p \left(1 - \frac{\ell(p^m)}{p^m}\right) &+ 
O(N^{2/(m+1)} (\log N)^C)
\\ 
&+ O(\# \{1\leq x\leq N : \exists p > N^{2/m+1} \text{ s.t. } p^m|Q(x)\}),
\end{aligned}\]
where
\[\ell(p^m) = \# \{x\in \mathbb{Z}/p^m : \mathfrak{p}^m|Q(x) 
\text{ for some $\mathfrak{p}\in I_K$ above $p$}\} .\]
The implied constant and $C$ depend only on $K$ and $Q$.
\end{prop}
\begin{proof}
We define a soil $(\mathscr{P},\mathscr{A},r,f)$ by
\[\begin{aligned}
\mathscr{P} &= \{\mathfrak{p}\in I_K : \text{$\mathfrak{p}$ prime}\},\;\;\;
\mathscr{A} = \{1,2,\dotsc,N\},\\
r(a) &= \{\mathfrak{p}\in P : \mathfrak{p}^m |Q(a)\},\;\;\;
f(a,d) = \begin{cases} 1 &\text{if $d=\emptyset$,}\\
0 &\text{otherwise.}\end{cases}
\end{aligned}\]
Let $h(d)$ be the positive integer generating the ideal
$(\prod_{\mathfrak{p}\in d} \mathfrak{p}^m)\cap \mathbb{Z}$.
Properties (h1) and (h2) are clear. 
Lemma \ref{lem:dontaskmemyadvice} gives us properties (A1) and (A2) with
\[\begin{aligned}\label{eq:agam}
X &= N,\;\;\; g(d_1,d_2) = 0 \;\text{ if $d_2$ non-empty,}\\
g(d,\emptyset) &= \# \{x\in \mathbb{Z}/h(d) : \mathfrak{p}^m | Q(x) \,
 \forall \mathfrak{p} \in d\},\\
|r_{d,\emptyset}| &\leq |\Disc Q|^3 (\deg Q)^{\# d},\;\;
|r_{d_1,d_2}| = 0 \;\text{ if $d_2$ non-empty,}\\
C_0 &\ll |\Disc Q|^3,\;\;\;
C_1, C_2 \ll \deg Q,\;\;\;
M_0 = N .
\end{aligned}\]
By Corollary \ref{cor:yugo},
$\sum_{a\in \mathscr{A}} f(a,r(a)) - N \sum_{d\subset \mathscr{P}}
\sum_{d'\subset d} \mu(d-d') \frac{g(d,d')}{h(d)}$
is at most
\[\begin{aligned}
N &\cdot \mathop{\sum_{d\subset \mathscr{P}}}_{h(d)>M}
 \frac{1}{h(d)} (2^{\# d} + C_0 C_1^{\# d} (3^{\# d} + 3)) 
+ \mathop{\sum_{d\subset \mathscr{P}}}_{M < h(d) \leq M^2} 
C_0 C_2^{\#d} (3^{\#d} + 3)\\
&+ \mathop{\sum_{d\subset \mathscr{P}}}_{h(d)\leq M} |\Disc Q|^3
(\deg Q)^{\# d} +
\mathop{\sum_{p\in \mathscr{P}}}_{h(\{p\})>M^2} S_{\{p\}} .
\end{aligned}\]
Given a positive integer $n$, there are at most 
$2^{(\deg K/\mathbb{Q}) \omega(n)}$ elements $d$ of $\mathscr{P}$ such
that $h(d) = m$. Moreover, for every prime $p$ not ramified in $K/\mathbb{Q}$,
$p|h(d)$ implies $p^m|h(d)$. Hence
\[\begin{aligned}
N \cdot \mathop{\sum_{d\subset \mathscr{P}}}_{h(d)>M} \frac{c^{\# d}}{h(d)}
 &\ll N \cdot \sum_{n>M^{1/m}} \frac{(2^{\deg K/\mathbb{Q}} c)^{\omega(n)}}{
n^m} \ll N \frac{(\log M)^{2^{2^{\deg K/\mathbb{Q}} c} - 1}}{N^{(m-1)/m}},\\
\mathop{\sum_{d\subset \mathscr{P}}}_{M<h(d)\leq M^2} c^{\# d} &\ll
 M^{2/m} (\log M)^{2^c - 1} ,\;\;\;
\mathop{\sum_{d\subset \mathscr{P}}}_{h(d)\leq M} c^{\# d} \ll
M^{1/m} (\log M)^{2^c - 1} .
\end{aligned}\]
We choose $M = N^{m/m+1}$, since then $\frac{N}{M^{(m-1)/m}} = M^{2/m}$.
It remains to show that
\[ \sum_{d\subset \mathscr{P}} \sum_{d'\subset d} \mu(d-d') \frac{g(d,d')}{
h(d)} 
= 
\prod_p \left(1 - \frac{\ell(p^m)}{p^m}\right) .\]
By (\ref{eq:agam}),
\[\begin{aligned}
\sum_{d\subset \mathscr{P}} &\sum_{d'\subset d} \mu(d-d') 
\frac{g(d,d')}{h(d)} = \sum_{d\subset \mathscr{P}} \mu(d)
\frac{g(d,0)}{h(d)} \\
&= \prod_{\text{$p$ prime}} \mathop{\sum_{d\subset \mathscr{P}}}_{h(d)=p}
\mu(d) \frac{g(d,0)}{h(d)} = \prod_{\text{$p$ prime}} \left(1 -
\frac{\ell(p^m)}{p^m}\right) \end{aligned}\]
and so we are done.
\end{proof}
While the following result could be stated in the same generality as
Prop. \ref{prop:owerr}, we restrict ourselves to the rationals 
and to square divisors in order
to avoid inessential complications.
\begin{prop}\label{prop:aogi}
Let $Q\in \mathbb{Z}\lbrack x,y\rbrack$ be a 
homogeneous polynomial. Then the number of
pairs of positive integers $x,y \leq N$ for which $Q(x,y)$ is free of
$m$th powers equals
\[\begin{aligned}
N^2 \prod_p \left(1 - \frac{\ell(p^2)}{p^4}\right) &+ 
O(N^{\frac{3}{2}} (\log N)^C) \\ &+
O(\{1\leq x,y\leq N : \exists p>N \text{ s.t. } p^2|Q(x,y)\}),
\end{aligned}\]
where $\ell(p^2) = \# \{ x,y \in \mathbb{Z}/p^2 : p^2 |Q(x,y)\}$.
The implied constant and $C$ depend only on $Q$.
\end{prop}
\begin{proof}
We define a soil $(\mathscr{P},\mathscr{A},r,f)$ by
\[\begin{aligned}
\mathscr{P} &= \{p\in \mathbb{Z}^+ : \text{$p$ prime}\},\;\;\;
\mathscr{A} = \{1\leq x,y\leq N : \gcd(x,y) = 1\},\\
r(a) &= \{p \in \mathscr{P} : p^2 | Q(a) \},\;\;
f(a,d) = \begin{cases} 1 &\text{if $d\ne \emptyset$,}\\
0 &\text{otherwise.} \end{cases}\end{aligned}\]
Let $h(d) = \prod_{p\in d} p^2$.
Properties (h1), (h2), (A1) and (A2) hold by 
Lemmas \ref{lem:sollat}, \ref{lem:penult} and \ref{lem:ramsay} with
\[\begin{aligned}
X &= N^2 \prod_{p} \left(1 - \frac{1}{p^2}\right),\;\;
g(d_1,d_2) = 0 \;\text{ if $d_2$ non-empty,}\\
g(d,\emptyset) &= 
\frac{\#\{x,y\in \mathbb{Z}/h(d): p^2|Q(x,y) \wedge 
(p\nmid x \vee p\nmid y) \; \forall p\in d\}}{
\prod_{p\in d} (p^2 - 1)},\\
|r_{d,\emptyset}| &\leq |\Disc Q|^3 (2 \deg Q)^{\# d} N \log N,\;\;
r_{d_1,d_2} = 0 \;\text{ if $d_2$ non-empty,}\\
C_0 &\ll |\Disc Q|^3,\;\;\;
C_1, C_2 = 2 \deg Q,\;\;\;
M_0 = N.
\end{aligned}\]
By Corollary \ref{cor:yugo},
$\sum_{a\in \mathscr{A}} f(a,r(a)) - X \sum_{d\in \mathscr{P}} 
\sum_{d' \subset d} \mu(d-d') g(d)$
is at most
\[\begin{aligned}
N^2 &\prod_p \left(1 - \frac{1}{p^2}\right) \sum_{n>M^{1/2}} \frac{1}{n^2} (2^{\omega(n)} +
 C_0 C_1^{\omega(n)} (3^{\omega(n)} + 3)) 
+
\mathop{\sum_{\text{$p$ prime}}}_{p>M} S_{\{p\}} 
\\
&+ \sum_{M^{1/2} < n \leq M} C_0 C_2^{\omega(n)} (3^{\omega(n)} + 3) +
\sum_{n\leq M^{1/2}} |\Disc Q|^3 (2 \deg Q)^{\omega(n)} N \log N ,\end{aligned}\]
which is in turn at most
\[N^2 \frac{(\log M)^C}{M^{1/2}} + M (\log M)^C + M^{1/2} N \log N +
\mathop{\sum_{\text{$p$ prime}}}_{p>M} S_{\{p\}} 
\] for some $C$ given by $C_0$, $C_1$, $C_2$. We set $M = N$.
Now
\[
\sum_{d\subset \mathscr{P}} \sum_{d'\subset d} \mu(d-d') \frac{g(d,d')}{h(d)}
 = \sum_{d\subset \mathscr{P}} \mu(d) \frac{g(d,\emptyset)}{h(d)} 
= \prod_p \left(1 - \frac{\ell'(p^2)}{p^{4}}\right) ,\]
where $\ell'(p^2) = \# \{x,y\in \mathbb{Z}/p^2 : p^2|Q(x,y) \wedge
(p\nmid x \vee p\nmid y)\}/(p^2-1)$. 
Clearly $\ell'(p^2) (p^2 - 1) = \ell(p^2) - p^2$. Hence
\[\begin{aligned} \prod_p \left(1 - \frac{\ell'(p^2)}{p^2}\right) &=
\prod_p \left(1 - \frac{\ell(p^2) - p^2}{p^2 (p^2 - 1)}\right) \\
&= \prod_p \frac{p^4 - \ell(p^2)}{p^2 (p^2 - 1)} = 
\prod_p \left(1 - \frac{\ell(p^2)}{p^4}\right) / 
\prod_p \left(1 - \frac{1}{p^2}\right) .\end{aligned}\]
The statement follows.
\end{proof}
\subsection{Sampling and averaging}\label{subs:samp}
A slight layer of abstraction is now called for. 
Let $J_0$ be an index set. Let $V_j$, $j\in J_0$, 
be measure spaces with positive measures $\mu_j$, $j\in J_0$. 
Consider a countable set $Z$ together with 
injections $\iota_j:Z\to V_j$ for every $j\in J_0$,
and a finite subset $Z_j \subset Z$ for every $n\in \mathbb{Z}^+$.
Let 
$\mathscr{M}_j$ be a collection
of measurable subsets of $V_j$, each of them of finite measure.
Let $\upsilon:\mathbb{Z}^+\to \mathbb{R}_0^+$
be a function with $\lim_{n\to \infty} \upsilon(n) = 0$. Assume that
\begin{equation}\label{eq:agty}
\left| \frac{1}{\# Z_n} \sum_{x\in Z_n \cap \bigcap_{j\in J} 
\iota_j^{-1}(M_j)}
1 - \prod_{j\in J} \mu_j(M_j) \right| \leq \upsilon(n) 
\end{equation}
for every finite subset $J\subset J_0$ and every tuple
$\{M_j\}_{j\in J}$, $M_j\in \mathscr{M}_j$. 

We may call $Z$ the {\em sample frame}, and 
$\mathscr{M}_j$, $j\in J_0$, the {\em sampling spaces}. A tuple
$(J_0,\{V_j\}, Z, \{\iota_j\}, 
\{Z_n\}, \upsilon, \{\mathscr{M}_j\})$
satisfying the conditions above, including (\ref{eq:agty}), 
 will be called a {\em sampling datum}.
Given a measure $\sigma_j$ on $V_j$ and a function $s_j:Z\to \mathbb{C}$
for every $j\in J$,
we say $(\{s_j\},\{\sigma_j\})$ 
is a {\em distribution pair} if $\max(|s_j(x)|)\leq 1$,
$\max |\sigma_j/\mu_j| \leq 1$, and
\[
\left| \frac{1}{\# Z_n} \sum_{x\in Z_n \cap \bigcap_{j\in J} 
\iota_j^{-1}(M_j)}
s(x) - \prod_{j\in J} \sigma_j(M_j) \right| \leq \upsilon_0(n)
\prod_{j\in J} \mu_j(M_j) + \upsilon_1(n) 
\]
for every finite subset $J\subset J_0$ and every tuple
$\{M_j\}_{j\in J}$, $M_j\in \mathscr{M}_j$, where 
$\upsilon_0(n),\upsilon_1(n):\mathbb{Z}^+\to \mathbb{R}_0^+$
are functions with $\lim_{n\to \infty} \upsilon_0(n) = 
\lim_{n\to \infty} \upsilon_1(n)=0$ and we write $s(x)$ for
$\prod_{j\in J} s_j(x)$. 

Given a positive integer $c$, 
we define $\mathscr{M}_{c,j}$ to be the collection of all sets
of the form \[(M_1\cup M_2 \cup \dotsb \cup M_{n_1}) - 
(M_1' \cup \dotsb \cup M_{n_2}'),\]
where $M_1,\dotsc,M_{n_1}\in \mathscr{M}_j$ and 
$M_1',\dotsc,M_{n_2}'\in \mathscr{M}_j$ and $n_1 + n_2 \leq c$.
It should be clear that (\ref{eq:agty}) and the inequality following it
hold for $M_j\in \mathscr{M}_{c,j}$ if the terms on the right are 
multiplied by $c^{\# J}$.

Let $j$ be an element of $J_0$. Let $M_j\in \mathscr{M}_j$.
Let $m:\mathbb{R}^+\to \mathbb{R}^+$ be a decreasing function with
$\int_0^{\infty} m(x)\, dx < \infty$. Let $c$ be a positive integer.
A function $f:M_j\to \mathbb{C}$
is {\em $(c,m)$-approximable} if there is a partition 
$M_j = M_{j,0} \cup M_{j,1} \cup M_{j,2} \cup \dotsb$ and a sequence
$\{y_i\}_{i\geq 1}$ of complex numbers such that
\begin{enumerate}
\item $M_{j,i} \in \mathscr{M}_{j,i}$ for $i\geq 1$,
\item $f(x) = y_i$ for $x\in U_i$, $i\geq 1$,
\item $\mu_j(M_{j,0}) = 0$,
\item $\mu_j(M_{j,i}) \leq m(i)$ for every $i\geq 1$.
\end{enumerate}

\begin{lem}\label{lem:aganel}
Let $(J_0,\{V_j\}, Z, \{\iota_j\}, \{Z_n\}, \upsilon, \{\mathscr{M}_j\})$ 
be a sampling datum, $(s_j,\sigma_j)$ a distribution pair.
Let $J$ be a finite subset of $J_0$. For every $j\in J$, choose
$M_j\in \mathscr{M}_j$ and
let $u_j:M_j\to \mathbb{C}$ be a $(c,m_j)$-approximable function with
$\max_x |u_j(x)| \leq 1$. Write $u(x) = \prod_{j\in J} u_j(x)$,
$s(x) = \prod_{j\in J} s_j(x)$, $c_0 = c^{\# J}$. Then,
for every subset $S$ of $(\mathbb{Z}^+)^J$,
\[\left| 
\frac{1}{\# Z_n} \sum_{x\in Z_n \cap \bigcap_{j\in J} 
\iota_j^{-1}(M_j)}
s(x) u(x)
- \prod_{j\in J} \int_{M_j} u_j(x) d\sigma_j\right|\]
is at most
\begin{equation}\label{eq:umbr}\begin{aligned}
\sum_{(t_j)\in (\mathbb{Z}^+)^J - S}
 \prod_{j\in J} m_j(t_j) &+
c_0 \upsilon_0(n) \prod_{j\in J} \mu_j(M_j) \\ &+
\# S \cdot c_0 \upsilon_1(n)
+ (\# S + 1) c_0 \upsilon(n).\end{aligned}
\end{equation}
\end{lem}
\begin{proof}
Since $|\sigma_j/\mu_j|\leq 1$, $|u_j|\leq 1$ and $u_j$ is $m_j$-approximable,
\[\begin{aligned}
&\left|\prod_{j\in J} \int_{M_j} u_j(x) d\sigma_j
- \sum_{(t_j)\in S} \prod_{j\in J} \int_{M_{j,t_j}} u_j(x) d\sigma_j
\right| \\&=
\left|\sum_{(t_j)\in (\mathbb{Z}^+)^J - S}
\prod_{j\in J} \int_{M_{j,t_j}} u_j(x) d\sigma_j \right|
\leq \sum_{(t_j)\in (\mathbb{Z}^+)^J - S}
 \prod_{j\in J} m_j(t_j),\end{aligned}\]
whereas
$\prod_{j\in J} \int_{M_{j,t_j}} u_j(x) d\sigma_j = \prod_{j\in J} y_j
 \cdot \sigma_j(M_{j,t_j})$. Now
\[ \left|\prod \sigma_j(M_{j,t_j}) -
\frac{1}{\# Z_n} \sum_{x\in Z_n \cap \bigcap_{j\in J} 
\iota_j^{-1}(M_{j,t_j})} s(x)\right|\]
is at most \[c_0 \upsilon_0(n)
\prod_{j\in J} \mu_j(M_{j,s_j}) + c_0\cdot \upsilon_1(n) . \]
Clearly
\[\begin{aligned}
\sum_{(t_j)\in S} c_0 \upsilon_0(n) \prod_{j\in J} \mu_j(M_{j,t_j}) &=
 c_0 \upsilon_0(n) \sum_{(t_j)\in (\mathbb{Z}^+)^J} \prod_{j\in J} 
 \mu_j(M_{j,t_j}) \\
&= c_0 \upsilon_0(n) \prod_{j\in J} \mu_j(M_j) .\end{aligned}\]
By (\ref{eq:agty}),
\[\begin{aligned}
\left|\sum_{(t_j)\notin S}
\frac{1}{\# Z_n} \sum_{x\in Z_n \cap \bigcap_{j\in J} 
\iota_j^{-1}(M_{j,t_j})} s(x) \right|
&\leq
\sum_{(t_j)\notin S}
 \frac{1}{\# Z_n} \sum_{x\in Z_n \cap \bigcap_{j\in J} 
\iota_j^{-1}(M_{j,t_j})} 1 \\
&=
\frac{\# (Z_n \cap \bigcap_{j\in J} \iota_j^{-1}(M_j))}{\# Z_n} \\ &-
\sum_{(t_j)\in S}
 \frac{1}{\# Z_n} \sum_{x\in Z_n \cap \bigcap_{j\in J} 
\iota_j^{-1}(M_{j,t_j})} 1 \\
&= \prod_{j\in J} \mu_j(M_j) -
\sum_{(t_j)\in S} \prod_{j\in J} \mu_j(M_{j,t_j}) \\ &+ 
O(c_0\cdot (\#S + 1) \upsilon(n)) \\
&= \sum_{(t_j)\in (\mathbb{Z}^+)^J - S} \prod_{j\in J} \mu_j(M_{j,t_j}) \\ & + 
O(c_0\cdot (\#S + 1) \upsilon(n))
.\end{aligned}\]
Both implied constants have absolute value at most $1$. The statement follows.
\end{proof}
It is time to produce some examples of sampling data.

\begin{Exa}\label{ex:adim}
Let $J_0$ be the set of rational primes. Let $V_p = \mathbb{Z}_p$,
$Z = \mathbb{Z}$, $Z_n = \{1,2,\dotsc,n\}$, 
$\iota_p:\mathbb{Z}\to \mathbb{Z}_p$ the natural injection,
$\mathscr{M}_p$ the set of additive subgroups of $\mathbb{Z}_p$ of
finite index, $\upsilon(n) = 1/n$. Then 
\[(J_0,\{V_p\},Z,\{\iota_p\},\{Z_n\},\upsilon,\{\mathscr{M}_p\})\] is
a sampling datum.
\end{Exa}

\begin{Exa}\label{ex:middim}
Let $J_0$ be the set of rational primes. Let $S\subset \mathbb{R}^2$
be a sector. Let
$V_p = \mathbb{Z}_p^2$,
 $Z = \mathbb{Z}^2 \cap S$,
$Z_n = \{-n\leq x,y\leq n\}\cap S$, $\iota_p:\mathbb{Z}\to \mathbb{Z}_p$
the natural injection,
$\upsilon(n) = O(1/n)$,
 $\mathscr{M}_p$ the set of all additive subgroups 
of $\mathbb{Z}_p^2$ of finite index.
Then
$(J_0,\{V_p\},Z,\{\iota_p\},\{Z_n\},\upsilon,\{\mathscr{M}_p\})$ is
a sampling datum. The implied constant in $\upsilon(n) = O(1/n)$ is absolute.
\end{Exa}

\begin{Exa}\label{ex:bdim}
Let $J_0$ be the set of rational primes. Let $S\subset \mathbb{R}^2$ be
a sector. Let
$V_p = \mathbb{Z}_p^2 - p \mathbb{Z}_p^2$, 
 $Z = \{x,y\in \mathbb{Z}: \gcd(x,y)=1\}\cap S$,
$Z_n = \{-n\leq x,y\leq n: \gcd(x,y)=1\}\cap S$,
$\upsilon(n) = O(\log n/n)$,
 $\mathscr{M}_p$ the collection of sets of the form $L\cap Z$, where by $L$
we mean an additive subgroup of $\mathbb{Z}_p^2$ of finite index.
By Lemma \ref{lem:penult}, 
$(J_0,\{V_p\},Z,\{\iota_p\},\{Z_n\},\upsilon,\{\mathscr{M}_p\})$ is
a sampling datum. The implied constant in 
$\upsilon(n) = O(\log n/n)$ is absolute.
\end{Exa}

The datum in example \ref{ex:bdim} is more natural than that in example
\ref{ex:middim} when the pairs $(x,y)\in \mathbb{Z}^2$ are meant
to represent rational numbers $x/y$.

\subsection{Averages}\label{subs:avmul}
We can now give explicit analogues of Proposition \ref{prop:unfair}. Our
conditions on the behavior of $u_p$ are fairly strict; somewhat laxer
conditions can be adopted with a consequent degradation in the
quality of the bounds. 

\begin{prop}\label{prop:poncho}
For every prime $p$,
let $u_p:\mathbb{Z}_p \to \mathbb{C}$ be a $(c_0,c_1 p^{-c_2 j})$-approximable
function, where
$c_0$ is a positive integer and $c_1$, $c_2$ are positive real numbers.
Assume 
that 
 $|u_p(x)|\leq 1$ for every $x\in \mathbb{Z}_p$. Assume, furthermore, that
$u_p(x) = 1$ unless $p^2|P(x)$, where $P\in \mathbb{Z}\lbrack x\rbrack$ is 
a fixed square-free polynomial. 
Then
\[\begin{aligned}
\frac{1}{N} \sum_{n=1}^N \prod_p u_p(n) &= \prod_p \int_{\mathbb{Z}_p}
 u_p(x) d x + O((\log N) N^{-1/3}) \\
&+ \frac{1}{N} O(\{1\leq x\leq N : \exists p > N^{2/3}
\text{ s.t. } p^2 | P(x)\}
),\end{aligned}\]
where the implied constant depends only on $P$, $c_0$, $c_1$ and $c_2$.
\end{prop}
\begin{proof}
We define a soil $(\mathscr{P},\mathscr{A},r,f)$ by
\[\begin{aligned}
\mathscr{P} &= \{p\in \mathbb{Z}^+ : \text{$p$ prime}\},\;\;\;
\mathscr{A} = \{1,2,\dotsc,N\},\\
r(a) &= \{p\in \mathscr{P} : p^2 | P(a)\},\;\;\;
f(a,d) = \prod_{p\in d} u_p(a) .\end{aligned}\]
Let $h(d) = \prod_{p\in d} p^2$. Properties (h1) and (h2) are clear.
Lemma \ref{lem:dontaskmemyadvice} 
gives us property (A1) with
\[X=N,\; C_0 \ll |\Disc P|^3,\; C_1,C_2\ll \deg P.\]
Choose the sampling datum in Example \ref{ex:adim} with $n=N$.
By Lemma \ref{lem:aganel}, property (A2) then holds with
\[\begin{aligned}
g(d_1,d_2) &= h(d_1) \cdot \prod_{p\in d_2} \int_{D_p} u_p(x) d \mu_p \cdot
\prod_{p\in d_1-d_2} \int_{D_p} d \mu_p \\
&= h(d_1) \cdot \prod_{p\in d_1 - d_2} \mu(D_p) \cdot \prod_{p\in d_2}
\int_{D_p} u_p(x) d \mu_p ,\end{aligned}\]
where $D_p = \{x\in \mathbb{Z}_p : p^2 | P(x)\}$, 
\[\begin{aligned}
r_{d_1,d_2} &\leq N\cdot \sum_{(t_p)\in (\mathbb{Z}^+)^{d_1} - S} \prod_{p \in 
d_1} m_p(t_p) + c_0^{\# d_1} \sum_{(t_p)\in S} 1 + c_0^{\# d_1} 
\cdot (\# S + 1) \\
&=
N\cdot \sum_{(t_p)\in (\mathbb{Z}^+)^{d_1} - S} \prod_{p \in 
d_1} m_p(t_p) + c_0^{\# d_1} \cdot (2 \# S + 1) \end{aligned}\]
for every $S \in (\mathbb{Z}^+)^{d_1}$ and $m_p(j) = c_1 p^{-c_2 j}$. (We can set $M_0$ arbitrarily large.)
Choose
\begin{equation}\label{eq:haber}
S = \vartimes_{p\in d_1} \{n\in \mathbb{Z}: 
1\leq n\leq \frac{c_2 \log N}{\log p}\} .\end{equation}
Then
\[\begin{aligned}
\sum_{(s_p) \in (\mathbb{Z}^+)^{d_1} - S} \prod_{p\in d_1} m_p(s_p) &\leq
 \sum_{p\in d_1} \sum_{j> \frac{c_2 \log N}{\log p}} m_p(j) \\
&\leq \sum_{p\in d_1} \frac{c_1}{1 - p^{-c_2}} 
e^{-\frac{c_2 \log N}{\log p} \cdot \frac{\log p}{c_2}} 
= \frac{c_1}{1 - 2^{-c_2}} \cdot \# d_1 / N ,\end{aligned}\]
and so
\[r_{d_1,d_2} \leq 2 c_0^{\# d_1} \prod_{p\in d_1} \left(
\frac{c_2 \log N}{\log p}\right) +
\frac{c_1}{1 - 2^{-c_2}} \cdot \# d_1 + c_0^{\# d_1} .\]
We can now apply Corollary \ref{cor:yugo}. We obtain that the absolute
value of the difference 
\[\begin{aligned}
\sum_{a\in \mathscr{A}} &f(a,r(a)) - N \cdot 
\sum_{d\in \mathscr{P}} \sum_{d' \subset d} \mu(d - d') \frac{g(d,d')}{h(d)}\\
&= \sum_{n=1}^N \prod_p u_p(n) - N \cdot
\sum_{\text{$d$ sq-free}} \sum_{d'|d} \mu(d/d') \prod_{p|d/d'} \mu_p(D_p)
 \prod_{p|d'} \int_{D_p} u_p(x) d \mu_p \\
&= \sum_{n=1}^N \prod_p u_p(n) - N \cdot \prod_p \int_{\mathbb{Z}_p}
 u_p(x) d\mu_p \end{aligned}\]
is at most (\ref{eq:vivsav}) with $C_3 = C_4 = 1$ and $M$ arbitrary.
The first term of (\ref{eq:vivsav}) is 
$O\left(\frac{(\log M)^c}{M^{1/2}} N\right)$, where $c$ and the implied 
constant depend only on $C_0$, $C_1$ and $C_2$. The second term is 
$O(M (\log M)^c)$. By Lemma \ref{lem:grogor}, the third term is no greater 
than
\[
\mathop{\sum_{1\leq d\leq M^{1/2}}}_{\text{$d$ sq-free}}
\sum_{d'|d} \left( 2 c_0^{\omega(d)} \prod_{p|d} \frac{c_2 \log N}{\log p} +
\frac{c_1}{1 - p^{-c_2}} \cdot 
\omega(d) + c_0^{\omega(d)} \right) \ll M^{1/2 + \epsilon} .\]
Set $M = N^{2/3}$. The result follows.
\end{proof}
\begin{prop}\label{prop:viento}
For every $p$, let 
$u_p : \mathbb{Z}_p^2 - p \mathbb{Z}_p^2 \to \mathbb{C}$ be a
$(c_0,c_1 p^{-c_2})$-approximable
function, where
$c_0$ is a positive integer and $c_1$, $c_2$ are positive real numbers.
Assume 
that 
 $|u_p(x,y)|\leq 1$ for all $(x,y)\in \mathbb{Z}_p^2 - p \mathbb{Z}_p^2$. 
Assume, furthermore, that
$u_p(x,y) = 1$ for all $(x,y)\in \mathbb{Z}_p^2 - p \mathbb{Z}_p^2$
such that
$p^2\nmid P(x)$, where $P\in \mathbb{Z}\lbrack x,y\rbrack$ is 
a fixed square-free homogeneous polynomial. 

Let $S\subset \mathbb{R}^2$ be a sector. Then
\[
\mathop{\mathop{\sum_{-N\leq x,y\leq N}}_{(x,y)\in S}}_{\gcd(x,y)=1}
\prod_p u_p(x,y)\] equals
 \[\begin{aligned} (&\# (S\cap \lbrack -N,N\rbrack^2)) \cdot
\prod_p \int_{\mathbb{Z}_p^2 - p \mathbb{Z}_p^2}
 u_p(x,y) d x d y + O((\log N) N^{4/3})\\ 
&+ O(\{x,y\in \lbrack -N,N\rbrack : 
 \gcd(x,y)=1, \exists p > N^{4/3}
\text{ s.t. } p^2 | P(x,y)\}),\end{aligned}\]
where the implied constants depend only on $P$, $c_0$, $c_1$ and $c_2$.
\end{prop}
\begin{proof}
We define a soil $(\mathscr{P},\mathscr{A},r,f)$ by
\[\begin{aligned}
\mathscr{P} &= \{p\in \mathbb{Z}^+ : \text{$p$ prime}\},\;\;\;
\mathscr{A} = \{1\leq x,y\leq N : \gcd(x,y)=1\},\\
r((x,y)) &= \{p\in \mathscr{P} : p^2 | P(x,y)\},\;\;\;
f((x,y),d) = \prod_{p\in d} u_p(x,y) .\end{aligned}\]
Let $h(d) = \prod_{p\in d} p^2$. Properties (h1) and (h2) are clear.
Lemmas \ref{lem:sollat} and \ref{lem:ramsay} give us property (A1) with 
\[X = N^2 \prod_p \left(1 - \frac{1}{p^2}\right),\;\;
C_0 \ll |\Disc P|^3,\;\; C_1,C_2 \ll \deg P .\]
Choose the sampling datum in Example \ref{ex:bdim} with $n=N$.
Proceed as in Prop. \ref{prop:poncho}.
\end{proof}
It is simple to show that certain natural classes of functions $u_p$
satisfy the conditions in Propositions \ref{prop:poncho} and 
\ref{prop:viento}.
\begin{lem}\label{lem:atra}
Let $P\in \mathbb{Z}\lbrack x\rbrack$ be a square-free polynomial.
Let $p$ be a rational prime. Let $u_p:\mathbb{Z}_p \to \mathbb{C}$
be such that $u_p(x)$ depends only on $x \mo p$ and $v_p(P(x))$.
Then $u_p$ is $(c_0, c_1 p^{-c_2 j})$-approximable, where
$c_0$, $c_1$ and $c_2$ depend only on $P$, not on $p$.
\end{lem}
\begin{proof}
Let $c_0 = 2 |\Disc P|^3 \deg P$,
$c_1 = |\Disc P|^3 \deg P$.
 Let $x_1,\dotsc,x_{c}$ be the solutions
to $P(x)\equiv 0 \mo p$ in $\mathbb{Z}/p$. Clearly $c\leq \deg P$.
Define \[M_{p,j_0,j_1} = \{x\in \mathbb{Z}_p : v_p(P(x)) = j_0,
x\equiv x_{j_1} \mo p\} .\]
Let $M_{p,0}=\{x\in \mathbb{Z}_p : p\nmid P(x)\}$. We have a partition
\[\mathbb{Z}_p = M_{p,0} \cup M_{p, 1,1} \cup \dotsb \cup M_{p,1,c} \cup
M_{p,2,1}\cup \dotsb \cup M_{p,2,c} \cup \dotsb .\]
By Lemma \ref{lem:dontaskmemyadvice}, 
\[\begin{aligned}
M_{p,j_0,j_1} &\in \mathscr{M}_{c_0,p},\\
\mu_p(M_{p,j_0,j_1}) &\leq \mu_p(\{x\in \mathbb{Z}_p : v_p(P(x)) = j_0\}) \leq
c_1 p^{-j_0} .\end{aligned}\]
Let $c_2 = 1/\deg P$. Thus $c_2\leq 1/c$.
The statement follows.
\end{proof}
\begin{lem}\label{lem:otro}
Let $P\in \mathbb{Z}\lbrack x,y\rbrack$ be a square-free homogeneous 
polynomial. Let $p$ be a rational prime.
Let $u_p:\mathbb{Z}_p^2 - p \mathbb{Z}_p^2 
\to \mathbb{C}$ be 
such that $u_p(x,y)$ depends only on $x y^{-1} \mo p \in \mathbb{P}^1(\mathbb{Z}/p \mathbb{Z})$ and $v_p(P(x,y))$.
Then $u_p$ is $(c_0, c_1 p^{-c_2 j})$-approximable, where
$c_0$, $c_1$ and $c_2$ depend only on $P$, not on $p$.
\end{lem}
\begin{proof}
Same as Lemma \ref{lem:atra}.
\end{proof}

\subsection{Averages with multipliers}\label{eq:avwis}

The framework developed in section \ref{subs:samp} involves two different sets
of measures, $\{\mu_j\}$ and $\{\sigma_j\}$, which have been set equal in the 
above. The fact that $\{\mu_j\}$ and $\{\sigma_j\}$ can be taken to be
different allows us to compute sums of the form $\sum s(n) \prod_p u_p(n)$,
where $s(n)$ is a function as in section \ref{subs:samp} -- a good 
{\em multiplier}, if you wish. 

The most natural non-trivial example may be $s(n) = \mu(n)$. The function
$\mu(n)$ averages to zero over arithmetic progressions. More precisely,
\begin{equation}\label{eq:aglok}
\frac{1}{N/m} \mathop{\sum_{1\leq n\leq N}}_{n\equiv a \mo m} \mu(n)
\ll N e^{- C (\log N)^{2/3} / (\log \log N)^{1/5}} \end{equation}
for $m\leq (\log N)^A$, where both
$C$ and the implied constant depend only on $A$.
We would like to show that $\mu(n) \prod_p u_p(n)$ averages to zero as well.
The utility of such a result follows from the discussion around
(\ref{eq:agort}): we want to average a function of the form
\[W(\mathscr{E}(n)) = \prod_p W_{0,p}(\mathscr{E}(n)) \cdot
\prod_p \frac{W_p(\mathscr{E}(n))}{W_{0,p}(\mathscr{E}(n))} =
f(n) \cdot \prod_p u_p(n) .\]
The factor $f(n)$ can generally be expressed in the form 
$\mu(P(n)) \cdot g(n)$, where $P$ is a polynomial (possibly constant)
and $g(n)$ is a factor that can be absorbed into the local factors
$u_p(n)$: 
\[g(n) \cdot \prod_p u_p(n) = \prod_p u_{p,1}(n) .\]
Thus, it befalls us to show that
$\mu(P(n)) \cdot \prod_p u_{p,1}(n)$
averages to zero. If $P$ is linear, we may use (\ref{eq:aglok}) and
Prop. \ref{prop:onkyo} below; in general, Prop. \ref{prop:onkyo} furnishes us
with a result conditional on 
\[\lim_{N\to \infty} \frac{1}{N} \mu(P(n)) = 0\;\;
 \text{        (Chowla's conjecture).}\] 

When we average $W_p(\mathscr{E}(t))$ over the rationals, we are computing
\[\mathop{\sum_{1\leq x,y\leq N}}_{\gcd(x,y)=1}
 W(\mathscr{E}(x/y)),\] or, in effect, the average of
\[\mu(P(x,y)) \prod_p u_{p,1}(x,y) .\]
For $\deg P = 1,2$, we have an analogue of (\ref{eq:aglok}). For 
$\deg P = 3$, we know from \cite{He2}, Th. 3.7.1 and Lem. 2.4.4, that
\[\mathop{\mathop{\sum_{-N\leq x,y\leq N}}_{(x,y)\in S \cap L}}_{\gcd(x,y)=1}
\mu(P(x,y)) \ll 
\frac{(\log \log N)^6 (\log \log \log N)}{\log N} \frac{N^2}{\lbrack
\mathbb{Z}^2 : L\rbrack} \]
for every lattice $L$ with $\lbrack \mathbb{Z}^2 : L \rbrack \leq (\log N)^A$,
where the implied constant depends only on $A$. 

A good multiplier need not have average zero. A simple and useful example of
a multiplier is given by the characteristic function of an arithmetic progression
$a + \mathbb{Z} m$ or of a lattice coset $L\subset \mathbb{Z}^2$.
\begin{prop}\label{prop:onkyo}
For every prime $p$,
let $u_p:\mathbb{Z}_p \to \mathbb{C}$ be a $(c_0,c_1 p^{-c_2 j})$-approximable
function, where
$c_0$ is a positive integer and $c_1$, $c_2$ are positive real numbers.
Assume 
that 
 $|u_p(x)|\leq 1$ for every $x\in \mathbb{Z}_p$. Assume, furthermore, that
$u_p(x) = 1$ unless $p^2|P(x)$, where $P\in \mathbb{Z}\lbrack x\rbrack$ is 
a fixed square-free polynomial. 

Let $s:\mathbb{Z}^+\to \mathbb{C}$ be a function with $|s(n)|\leq 1$ for
every $n\in \mathbb{Z}^+$. Let $\sigma_p$ be a measure on $\mathbb{Z}_p$
with $\max_{S\subset \mathbb{Z}_p} |\sigma_p(S)/\mu_p(S)| \leq 1$,
where $\mu_p$ is the usual measure on $\mathbb{Z}_p$.
Suppose there are $\epsilon(N)$, $\eta(N)$ with
 $0\leq \epsilon(N) \leq 1$,
$1\leq \eta(N) \leq N$, such that
\begin{equation}\label{eq:aglan}
\mathop{\sum_{1\leq n\leq N}}_{n\equiv a \mo m} s(n) =
N\cdot \prod_p \sigma_p(\{n\in \mathbb{Z}_p : n\equiv a \mo p^{v_p(m)}\})
+ O\left(\frac{\epsilon(N) N}{m}\right)\end{equation}
for all $a$, $m$ with $0<m\leq \eta(N)$. Then
\begin{equation}\label{eq:ogo}\begin{aligned}
\frac{1}{N} \sum_{n=1}^N s(n) \prod_p u_p(n) &=
\prod_p \int_{\mathbb{Z}_p} u_p(x) d \sigma_p
+ O\left(\epsilon(N) + \frac{1}{\eta(N)^{1/2-\epsilon'}}\right) \\
&+  \frac{1}{N} O(\{1\leq x\leq N: \exists p>N^{1/2} \text{\,s.t.\,} p^2|P(x)\}).
\end{aligned}\end{equation}
for all $\epsilon'>0$,
where the implied constant depends only on $P$, $c_0$, $c_1$, $c_2$,
$\epsilon'$,
and the implied constant in \ref{eq:aglan}.
\end{prop}
\begin{proof}
Define a soil $(\mathscr{P},\mathscr{A},r,f)$ as in Prop. \ref{prop:poncho}.
Choose the sampling datum in Example \ref{ex:adim} with $n=N$. Then (A2)
holds with $g(d_1,d_2)$ as in Prop. \ref{prop:poncho} with $\sigma_p$ instead
of $\mu_p$, and
\[\begin{aligned}
r_{d_1,d_2} &\leq N\cdot \sum_{(s_j)\in (\mathbb{Z}^+)^{d_1} - S} \prod_{p \in 
d_1} m_p(s_p) + c_0^{\# d_1} \epsilon(N) N \prod_{j\in d_1} \mu_j(M_j) \\ &+
\# S \frac{c_0^{\# d_1} N}{\eta(N)} + 
 (\# S + 1) c_0^{\# d_1} \end{aligned}\]
for every $S \in (\mathbb{Z}^+)^{d_1}$ and $m_p(j) = c_1 p^{-c_2 j}$. Choose
\[S = \vartimes_{p\in d_1} \{n\in \mathbb{Z} : 1\leq n\leq \frac{c_2 
\log \eta(N)}{\log p}\} .\]
Then
\[\sum_{(s_j)\in (\mathbb{Z}^+)^{d_1} - S} \prod_{p\in d_1} m_p(s_p)
\leq \frac{c_1}{1 - p^{-c_2}} \cdot \# d_1 /\eta(N) \]
and
\[r_{d_1,d_2} \ll c_0^{\# d_1} \cdot \left(
\frac{\epsilon(N) N}{h(d_1)} (\deg P)^{\# d_1} +
\frac{N}{\eta(N)} \# d_1 + \frac{N}{\eta(N)} \prod_{p\in d_1} 
\frac{c_2 \log \eta(N)}{\log p} \right) .\]
Apply Corollary \ref{cor:yugo}. The first and second terms of 
(\ref{eq:vivsav}) are as in Proposition \ref{prop:poncho}. By Lemma
\ref{lem:crudel}, the third term is
\[O\left(\epsilon(N) N + \frac{N M^{1/2 + \epsilon'}}{\eta(N)^{1-\epsilon'}}\right) .\]
The fourth term is
\[\begin{aligned}
\sum_{p>M} \# \{1\leq x\leq N: p^2|P(x)\} &=
 \sum_{M<p\leq N^{1/2}} \# \{1\leq x\leq N : p^2|P(x)\} \\ &+
 \sum_{p>N^{1/2}} \# \{1\leq x\leq N : p^2|P(x)\} \\
 &\ll \sum_{M<p\leq N^{1/2}} N/p^2 \\ &+
 \# \{1\leq x\leq N: \exists p>N^{1/2} \text{\,s.t.\,} p^2|P(x)\}\\
&\leq N M^{-1} \\ &+  \# \{1\leq x\leq N: \exists p>N^{1/2} \text{\,s.t.\,} p^2|P(x)\}.
\end{aligned}\]
Set $M = \eta(N)$. 
The result follows.
\end{proof}
\begin{prop}\label{prop:jvc}
For every $p$, let 
$u_p : \mathbb{Z}_p^2 - p \mathbb{Z}_p^2 \to \mathbb{C}$ be a
$(c_0,c_1 p^{-c_2})$-approximable
function, where
$c_0$ is a positive integer and $c_1$, $c_2$ are positive real numbers.
Assume 
that 
 $|u_p(x,y)|\leq 1$ for all $(x,y)\in \mathbb{Z}_p^2 - p \mathbb{Z}_p^2$. 
Assume, furthermore, that
$u_p(x,y) = 1$ for all $(x,y)\in \mathbb{Z}_p^2 - p \mathbb{Z}_p^2$
such that
$p^2\nmid P(x)$, where $P\in \mathbb{Z}\lbrack x,y\rbrack$ is 
a fixed square-free homogeneous polynomial. 

Let $S\subset \mathbb{R}^2$ be a sector. Let 
$s:\{(x,y)\in \mathbb{Z}^2 : \gcd(x,y)=1\} \to \mathbb{C}$ be a function
with $|s(x,y)|\leq 1$ for all $x,y\in \mathbb{Z}$, $\gcd(x,y)=1$.
Let $\sigma_p$ be a measure on $\mathbb{Z}_p^2 - p \mathbb{Z}_p^2$
with $\max_{S\subset \mathbb{Z}_p^2 - p \mathbb{Z}_p^2} 
|\sigma_p(S)/\mu_p(S)| \leq 1$,
where $\mu_p$ is the usual measure on $\mathbb{Z}_p^2$.
 Suppose
there are 
$\epsilon(N)$, $\eta(N)$ with
 $0\leq \epsilon(N) \leq 1$,
$1\leq \eta(N) \leq N$, such that
\begin{equation}\label{eq:oglorsk}
\mathop{\mathop{\sum_{-N\leq x,y\leq N}}_{(x,y)\in S\cap L}}_{\gcd(x,y)=1}
 s(x,y) = N^2 \prod_p 
\sigma_p(\{(x,y)\in (\mathbb{Z}_p^2 - p \mathbb{Z}_p^2) \cap L_p\})
+ O\left(\frac{\epsilon(N) N}{\lbrack \mathbb{Z}^2 : L \rbrack}\right)
\end{equation}
for all lattices $L\subset \mathbb{Z}^2$ with 
$\lbrack \mathbb{Z}^2 : L\rbrack \ll \eta(N)$, where $L_p\subset \mathbb{Z}_p^2$
is the additive subgroup
 of $\mathbb{Z}_p^2$ generated by $\mathbb{Z}_p L$. Then
\begin{equation}\label{eq:whocar}\begin{aligned}
\frac{1}{N^2} 
&\mathop{\mathop{\sum_{-N\leq x,y\leq N}}_{(x,y)\in S}}_{\gcd(x,y)=1}
s(x,y) \prod_p u_p(x,y) = \int_{\mathbb{Z}_p^2 - p \mathbb{Z}_p^2} u_p(x,y)
d \sigma_p\\
&+ O\left(\epsilon(N) + \frac{1}{\eta(N)^{1/2-\epsilon'}}\right) \\
&+ \frac{1}{N^2} O(\{-N\leq x,y\leq N : \gcd(x,y)=1, \exists p > N
\text{ s.t. } p^2 | P(x,y)\}) 
\end{aligned}\end{equation}
for all $\epsilon'>0$,
where the implied constant depends only on $P$, $c_0$, $c_1$, $c_2$,
$\epsilon'$, and
the implied constant in (\ref{eq:oglorsk}).
\end{prop}
\begin{proof}
As in Proposition \ref{prop:onkyo}. Use Prop. \ref{prop:viento} instead
of \ref{prop:poncho}. 
\end{proof}
The reader may wonder why Proposition \ref{prop:onkyo} requires information
on averages on arithmetic progressions (\ref{eq:aglan}), yet seems
to furnish data only on averages over all positive integers (\ref{eq:ogo}).
In fact, we can obtain information on averages over arithmetic progressions
by applying Proposition \ref{prop:onkyo} to a new multiplier $s_0$
defined in terms of a given arithmetic progression $a + m \mathbb{Z}$:
\[s_0(n) = \begin{cases} s(n) &\text{if $n\equiv a \mo m$,}\\
0 &\text{otherwise.}\end{cases}\]
The same can be done as far as Proposition \ref{prop:jvc} and equations
(\ref{eq:oglorsk}) and (\ref{eq:whocar}) are concerned, with lattice
and lattice cosets playing the role of arithmetic progressions.
\section{Large square divisors of values of polynomials}\label{sec:gust}
It is now time to estimate the error terms denoted by $\delta(N)$ in
the introduction. These are error terms of the form
\[\frac{1}{N} O(\{1\leq x\leq N : \exists p > N^{1/2}
\text{ s.t. } p^2 | P(x)\}\]
and
\[\frac{1}{N^2} O(\{-N\leq x,y\leq N : \gcd(x,y)=1, \exists p > N
\text{ s.t. } p^2 | P(x,y)\}),\]
where $P \in \mathbb{Z}\lbrack x \rbrack$ 
 (resp. $P\in \mathbb{Z}\lbrack x,y\rbrack$ homogeneous) is a given
polynomial square-free as an element of $\mathbb{Q}\lbrack x\rbrack$
(resp. square-free as an element of $\mathbb{Q}\lbrack x,y\rbrack$).
In order to go beyond previous estimates, we will need to go beyond
sieve theory into diophantine geometry.
\subsection{Elliptic curves, heights and lattices}\label{subs:shlime}
As is usual, we write $\hat{h}$ for the canonical height on an elliptic
curve $E$, and $h_x$, $h_y$ for the height on $E$ with respect to $x$, $y$:
\[h_x((x,y)) = \begin{cases} 0 &\text{if $P=O$,}\\
\log H(x) &\text{if $P=(x,y)$,}\end{cases}\]
\[h_y((x,y)) = \begin{cases} 0 &\text{if $P=O$,}\\
\log H(y) &\text{if $P=(x,y)$,}\end{cases}\]
where $O$ is the origin of $E$, taken to be the point at infinity, and
\[\begin{aligned}
H(y) &= (H_K(y))^{1/\lbrack K : \mathbb{Q}\rbrack},\\
H_K(y) &= \prod_v \max(|y|_v^{n_v},1),\end{aligned}\]
where $K$ is any number field containing $y$,  the product
$\prod_v$ is taken over all places $v$ of $K$, and $n_v$ denotes
the degree of $K_v/\mathbb{Q}_v$.

In particular, if $x$ is a rational number $x_0/x_1$, $\gcd(x_0,x_1)=1$, then
\[\begin{aligned}
H(x) &= H_{\mathbb{Q}}(x) = \max(|x_0|,|x_1|),\\
h_x((x,y)) &= \log(\max(|x_0|,|x_1|)).\end{aligned}\]

The differences $|\hat{h} - \frac{1}{2} h_x|$ and 
$|\hat{h} - \frac{1}{3} h_y|$ are bounded
on the set of all points of $E$ (not merely on $E(\mathbb{Q})$). This
basic property of the canonical height will be crucial in our analysis.

\begin{lem}\label{lem:soide}
Let $f\in \mathbb{Z}\lbrack x\rbrack$ be a cubic polynomial of 
non-zero discriminant. For every square-free rational integer $d$, 
let $E_d$ be
the elliptic curve 
\[E_d : d y^2 = f(x) .\]
Let $P = (x,y) \in E_d(\mathbb{Q})$. Consider the point
 $P'=(x,d^{1/2} y)$ on $E_1$.
Then $\hat{h}(P) = \hat{h}(P')$, where the canonical heights are defined on
$E_d$ and $E_1$, respectively,
\end{lem}
\begin{proof} Clearly $h_x(P') = h_x(P)$. Moreover $(P+P)' = P' + P'$.
Hence \[\hat{h}(P) = \frac{1}{2} \lim_{N\to \infty} 4^{-N}
h_x(\lbrack 2^N\rbrack P) = \frac{1}{2} \lim_{N\to \infty} 4^{-N}
h_x(\lbrack 2^N\rbrack P') = \hat{h}(P') .\]
\end{proof}

\begin{lem}\label{lem:pseusi}
Let $f\in \mathbb{Z}\lbrack x\rbrack$ be an irreducible polynomial.
 Let $C$ be the curve given by
$C:y^2 = f(x)$. Let $d\in \mathbb{Z}$ be square-free. Let $x$, $y$ be rational
numbers, $y\ne 0$, such that $P = (x, d^{1/2} y)$ lies on $E$. Then
\[h_y(P) = \log H(d^{1/2} y) \leq \frac{3}{8} \log |d| + C_f,\]
where $C_f$ is a constant depending only on $f$.
\end{lem}
\begin{proof}
Write $y = y_0/y_1$, where $y_0$ and $y_1$ are coprime integers. Then
\begin{equation}\label{eq:wraw}
H(y) = \max\left(\frac{|y_0| |d|^{1/2}}{\sqrt{\gcd(d,y_1^2)}},
  \frac{|y_1|}{\sqrt{\gcd(d,y_1^2)}}\right) .\end{equation}

Write $a$ for the leading coefficient of $f$. Let $p|\gcd(d,y_1^2)$, $p\nmid a$.
Since $d$ is square-free, $p^2\nmid \gcd(d,y^2)$. Suppose $p^2\nmid y_1$. Then $\nu_p(d y^2) = -1$.
However, $d y^2 = f(x)$ implies that, if $\nu_p(x)\geq 0$, then $\nu_p(d y^2)\geq 0$, and
if $\nu_p(x)<0$, then $\nu_p(d y^2)\leq -3$. Contradiction. Hence $p|\gcd(d,y_1^2)$,
$p\nmid a$ imply $p^2\nmid \gcd(d,y_1^2)$, $p^2|y_1$. Therefore
 $|y_1|\geq (\gcd(d,y_1^2)/a)^2$. 

By (\ref{eq:wraw}) it follows that
\[\begin{aligned}
H(P)&\geq \max\left(\frac{|d|^{1/2}}{\sqrt{\gcd(d,y_1^2)}},
  \frac{|y_1|}{\sqrt{\gcd(d,y_1^2)}}\right)\\
&\geq \max\left(\frac{|d|^{1/2}}{\sqrt{\gcd(d,y_1^2)}},
  \frac{(\gcd(d,y_1^2))^{3/2}}{a^2}\right) .\end{aligned}\]
Since $\max(|d|^{1/2} z^{-1/2}, z^{3/2}/a^2)$ is minimal when
$|d|^{1/2} z^{-1/2} = z^{3/2}/a^2$, i.e., when
$z = a |d|^{1/4}$, we obtain
\[H(P)\geq |d|^{3/8} |a|^{-1/2} .\]
Hence 
\[h_y(P) = \log H(P) \geq \frac{3}{8} \log |d| - \frac{1}{2} \log |a| .\]
\end{proof}
\begin{cor}\label{cor:nosilv}
Let $f\in \mathbb{Z}\lbrack x\rbrack$ be a cubic polynomial of non-zero discriminant.
For every square-free rational integer $d$, let $E_d$ be the elliptic curve
\[E_d : d y^2 = f(x) .\]
Let $P = (x,y)\in E_d(\mathbb{Q})$. Then
\[\hat{h}(P)\geq \frac{1}{8} \log |d| + C_f ,\]
where $C_f$ is a constant depending only on $f$.
\end{cor}
\begin{proof}
Let $P' = (x, d^{1/2} y)\in E_1$. By Lemma \ref{lem:soide}, 
$\hat{h}(P) = \hat{h}(P')$. The difference $|\hat{h} - h_x|$ is bounded on $E$. The
statement follows from Lemma \ref{lem:pseusi}.
\end{proof}

The following crude estimate will suffice for some of our purposes.
\begin{lem}\label{lem:salom}
Let $Q$ be a positive definite quadratic form on $\mathbb{Z}^r$. Suppose
$Q(\vec{x})\geq c_1$ for all non-zero $\vec{x}\in \mathbb{Z}^r$. Then there are
at most 
\[(1 + 2 \sqrt{c_2/c_1})^{r}\]
values of $\vec{x}$ for which $Q(\vec{x})\leq c_2$.
\end{lem}
\begin{proof}
There is a linear bijection $f:\mathbb{Q}^r\to \mathbb{Q}^r$ taking $Q$
to the square root of the Euclidean norm: $Q(\vec{x}) = |f(\vec{x})|^2$ for
all $\vec{x}\in \mathbb{Q}^r$. Because $Q(\vec{x})> c_1$ for all
non-zero $\vec{x}\in \mathbb{Z}^r$, we have that $f(\mathbb{Z}^r)$ is a lattice
$L\subset \mathbb{Q}^r$ such that $|\vec{x}|\geq c_1^{1/2}$ for all
$\vec{x}\in L$, $\vec{x}\ne 0$. We can draw a sphere $S_{\vec{x}}$ of
radius $\frac{1}{2} c_1^{1/2}$ around each point $\vec{x}$ of $L$. The
spheres do not overlap. If $\vec{x}\in L$, $|\vec{x}|\in c_2^{1/2}$,
then $S_{\vec{x}}$ is contained in the sphere $S'$ of radius 
$c_2^{1/2} + c_1^{1/2}/2$ around the origin. The total volume
of all spheres $S_{\vec{x}}$ within $S'$ is no greater than the volume of
$S'$. Hence
\[\# \{\vec{x}\in L: |\vec{x} | \leq c_2^{1/2}\} \cdot (c_1^{1/2}/2)^r
\leq (c_2^{1/2} + c_1^{1/2}/2)^r .\]
The statement follows.
\end{proof}
\begin{cor}\label{cor:magd}
  Let $E$ be an elliptic curve over $\mathbb{Q}$. Suppose there are no non-torsion points
$P\in E(\mathbb{Q})$ of canonical height $\hat{h}(P)<c_1$. Then there are at most
\[O\left((1 + 2\sqrt{c_2/c_1})^{\rnk(E)}\right)\]
points $P\in E(\mathbb{Q})$ for which $\hat{h}(P) < c_2$. The implied constant is absolute.
\end{cor} \begin{proof}
The canonical height $\hat{h}$ is a positive definite quadratic form on the free part
$\mathbb{Z}^{\rnk(E)}$ of $E(\mathbb{Q}) \sim \mathbb{Z}^{\rnk(E)} \times T$. A classical
theorem of Mazur's \cite{Maz} states that the cardinality of $T$ is at most $16$. 
 Apply Lemma \ref{lem:salom}.
\end{proof}
Note that we could avoid the use of Mazur's theorem, since Lemmas \ref{lem:soide} and 
\ref{lem:pseusi} imply that the torsion group of $E_d$ is either $\mathbb{Z}/2$ or trivial
for large enough $d$.

\subsection{Twists of cubics and quartics}

Let $f(x) = a_4 x^4 + a_3 x^3 + a_2 x^2 + a_1 x + a_0 \in 
\mathbb{Z}\lbrack x\rbrack$ be an irreducible polynomial of degree $4$.
For every square-free $d\in \mathbb{Z}$, consider the curve
\begin{equation}\label{eq:conco1} C_d : d y^2 = f(x) .\end{equation}
If there is a rational point $(r,s)$ on $C_d$, then there is a birational
map from $C_d$ to the elliptic curve
\begin{equation}
E_d : d y^2 = x^3 + a_2 x^2 + (a_1 a_3 - 4 a_0 a_4) x - 
(4 a_0 a_2 a_4 - a_1^2 a_4 - a_0 a_3^2) .\end{equation}
Moreover, we can construct such a birational map in terms of $(r,s)$
as follows. Let $(x,y)$ be a rational point on $C_d$. We can rewrite
(\ref{eq:conco1}) as
\[y^2 = \frac{1}{d} f(x) .\]
We change variables:
\[x_1 = x-r, \;\; y_1 = y\]
satisfy
\[y^2 = \frac{1}{d} \left(\frac{1}{4!} f^{(4)}(r) x_1^4 +
\frac{1}{3!} f^{(3)}(r) x_1^3 + \frac{1}{2!} f''(r) x_1^2 +
\frac{1}{1!} f'(r) x_1 + f(r)\right) .\]
We now apply the standard map for putting quartics in Weierstrass form:
\[\begin{aligned}
x_2 &= (2 s (y_1 + s) + f'(r) x_1/d)/x_1^2,\\
y_2 &= (4 s^2 (y_1 + s) + 2 s (f'(r) x_1/d + f''(r) x_1^2/(2 d)) -
(f'(r)/d)^2 x_1^2/(2 s))/x_1^3 \end{aligned}\]
satisfy
\begin{equation}\label{eq:pfauen1}
y_2^2 + A_1 x_2 y_2 + A_3 y_2 = x_2^3 + A_2 x_2^2 + A_4 x_2 + A_6
\end{equation}
with
\[\begin{aligned}
A_1 &= \frac{1}{d} f'(r)/s,\;\;&A_2 = 
                    \frac{1}{d}(f''(r)/2 - (f'(r))^2/(4 f(r))),\\
A_3 &= \frac{2 s}{d} f^{(3)}(r)/3!,\;\;&A_4 = -\frac{1}{d^2} \cdot
4 f(r) \cdot \frac{1}{4!} f^{(4)}(r) ,\\
A_6 &= A_2 A_4 .\end{aligned}\]
To take (\ref{eq:pfauen1}) to $E_d$, we apply a linear change of variables:
\[
x_3 = d x_2 + r (a_3 + 2 a_4 r),\;\;\;
y_2 = \frac{d}{2} (2 y_2 + a_1 x_2 + a_3)\]
satisfy
\[d y_3^2 = 
x_3^3 + a_2 x_3^2 + (a_1 a_3 - 4 a_0 a_4) x_3 - 
(4 a_0 a_2 a_4 - a_1^2 a_4 - a_0 a_3^2) .\]
We have constructed a birational map $\phi_{r,s}(x,y)\mapsto (x_3,y_3)$
from $C_d$ to $E_d$.

Now consider the equation
\begin{equation}\label{eq:threest}
d y^2 = a_4 x^4 + a_3 x^3 z + a_2 x^2 z^2 + a_1 x z^3 + a_0 z^4 .
\end{equation}
Suppose there is
a solution $(x_0,y_0,z_0)$ to (\ref{eq:threest}) with 
$x_0, y_0, z_0 \in \mathbb{Z}$, $|x_0|, |z_0|\leq N$, $z_0\ne 0$. Then 
$(x_0/z_0,y_0/z_0^2)$ is a rational point on (\ref{eq:conco1}). We can
set $r = x_0/z_0$, $s = y_0/z_0^2$ and define a map $\phi_{r,s}$ from
$C_d$ to $E_d$ as above. Now let $x, y, z\in \mathbb{Z}$, $|x|, |z|\leq N$,
$z_0\ne 0$,
be another solution to (\ref{eq:threest}). Then
\[P = \phi_{r,s}(x_0/z_0,y_0/z_0^2)\]
is a rational point on $E_d$. Notice that $|y_0|, |y| \ll (N^4/d)^{1/2}$.
Write \[\phi_{r,s}(P) = (u_0/u_1,v),\]
where $u_0,u_1\in \mathbb{Z}$, $v\in \mathbb{Q}$, $\gcd(u_0,u_1)=1$. By a 
simple examination of the construction of $\phi_{r,s}$ we can determine
that $\max(u_0,u_1) \ll N^7$, where the implied constant depends only
on $a_0, a_1, \dotsb, a_4$.
In other words,
\begin{equation}\label{eq:quadsta}
h_x(P) \leq 7 \log N + C,\end{equation}
where $C$ is a constant depending only on $a_j$. Notice that (\ref{eq:quadsta})
holds even for $(x,y,z) = (x_0,y_0,z_0)$, as then $P$ is the origin of $E$.

The value of $h_x(P)$ is independent of whether $P$ is considered as a rational
point of $E_d$ or as a point of $E_1$. Let $\hat{h}_{E_1}(P)$ be the 
canonical height of $P$ as a point of $E_1$. Then
\[|\hat{h}_{E_1} (P) - \frac{1}{2} h_x(P)|\leq C',\]
where $C'$ depends only on $f$. By Lemma \ref{lem:soide}, the canonical height
$\hat{h}_{E_1}(P)$ of $P$  as a point of $E_1$ equals the canonical height
$\hat{h}_{E_d}(P)$ of $P$ as a point of $E_d$. Hence 
\[|\hat{h}_{E_d}(P) - \frac{1}{2} h_x(P)|\leq C' .\]
Then, by (\ref{eq:quadsta}),\[\hat{h}_{E_d}(P)\leq \frac{7}{2} \log N + (C/2 + C') .\]

We have proven
\begin{lem}\label{lem:ari} Let 
$f(x,z) = 
a_4 x^4 + a_3 x^3 z + a_2 x^2 z^2 + a_1 x z^3 + a_0 z^4\in \mathbb{Z}\lbrack x,z\rbrack$ 
be an irreducible homogeneous polynomial. Then there is a constant $C_f$
such that the following holds. Let $N$ be any positive integer. Let $d$
be any square-free integer. Let $S_{d,1}$ be the set of all
solutions $(x,y,z)\in \mathbb{Z}^3$ to
\[d y^2 = f(x,z) \]
satisfying $|x|, |z| \leq N$, $\gcd(x,z)=1$. 
Let $S_{d,2}$ be the set of all rational points $P$
on
\begin{equation}\label{eq:vova} E_d : d y^2 = 
x^3 + a_2 x^2 + (a_1 a_3 - 4 a_0 a_4) x - 
(4 a_0 a_2 a_4 - a_1^2 a_4 - a_0 a_3^2) \end{equation}
with canonical height
\[\hat{h}(P) \leq \frac{7}{2}\log N + C_f .\]
Then there is an injective map from $S_{d,1}$ to $S_{d,2}$.
\end{lem}
We can now apply the results of subsection \ref{subs:shlime}.
\begin{prop}\label{prop:clinton}
Let $f(x,z) =
a_4 x^4 + a_3 x^3 z + a_2 x^2 z^2 + a_1 x z^3 + a_0 z^4\in 
\mathbb{Z}\lbrack x,z\rbrack$ 
be an irreducible homogeneous polynomial. Then there are constants 
$C_{f,1}$, $C_{f,2}$, $C_{f,3}$
such that the following holds. Let $N$ be any positive integer. Let $d$
be any square-free integer. Let $S_{d}$ be the set of all 
solutions $(x,y,z)\in \mathbb{Z}^3$ to
\[d y^2 = f(x,z) \]
satisfying $|x|, |z| \leq N$, $\gcd(x,z)=1$. Then
\[\# S_d \ll \begin{cases} 
\left(1 + 2 \sqrt{(\frac{7}{2} \log N + C_{f,1})/
(\frac{1}{8} \log |d| + C_{f,2})}\right)^{\rnk(E_d)}
&\text{if $|d|\geq C_{f,4}$,}\\
\left(1 + 2 C_{f,3} \sqrt{\frac{7}{2} \log N + C_{f,1}}\right)^{\rnk(E_d)}
&\text{if $|d| < C_{f,4}$,}\end{cases}\]
where $C_{f,4} = e^{9 C_{f,2}}$, $E_d$ is as in (\ref{eq:vova}),
and the implied constant depends only on $f$.
\end{prop}
\begin{proof}
If $|d|\leq C_{f,4}$, apply Corollary \ref{cor:magd} and Lemma
\ref{lem:ari}. If $|d|> C_{f,4}$, apply Corollary \ref{cor:nosilv},
Corollary \ref{cor:magd} and Lemma
\ref{lem:ari}. 
\end{proof}
\subsection{Divisor functions and their averages}
As is usual, we denote by $\omega(d)$ the number of prime divisors of
a positive integer $d$.
Given an extension $K/\mathbb{Q}$,
we define \[\omega_K(d) = 
\mathop{\sum_{\mathfrak{p}\in I_K}}_{\mathfrak{p}|d} 1 .\]
\begin{lem}\label{lem:notse}
Let $f(x)\in \mathbb{Z}\lbrack x\rbrack$ be an irreducible polynomial of
degree $3$ and non-zero discriminant.
Let $K = \mathbb{Q}(\alpha)$, where $\alpha$ is a root of $f(x)=0$.
For every square-free rational integer $d$, let $E_d$ be the elliptic curve given by
\[ d y^2 = f(x) .\]
Then 
\[\rnk(E_d) = C_f + \omega_K(d) - \omega(d),\]
where $C_f$ is a constant depending only on $f$.
\end{lem}
\begin{proof}
Write $f(x) = a_3 x^3 + a_2 x^2 + a_1 x + a_0$.
Let $f_d(x) = a_3 x^3 + d a_2 x^2 + d^2 a_1 x + d^3 a_0$. 
Then $d \alpha$ is a root of $f_d(x)=0$.
Clearly $\mathbb{Q}(d \alpha) = \mathbb{Q}(\alpha)$. If $p$ is
a prime of good reduction for $E_1$, then $E_d$ will have additive reduction
at $p$ if $p|d$, and good reduction at $p$ if $p\nmid d$.
 The statement now follows immediately from the standard bound in
 \cite{BK}, Prop. 7.1.
\end{proof}

If $d y^2 = f(x)$ is to have any integer points $(x,y)\in \mathbb{Z}^2$
at all, no prime unsplit in $\mathbb{Q}(\alpha)/\mathbb{Q}$ can divide
$d$. We define
\begin{equation}\label{eq:manandur}
R(\alpha,d) = \begin{cases} 2^{\alpha \omega_{\mathbb{Q}(\alpha)}(d) - 
\alpha \omega(d)} &\text{if no $p|d$ is unsplit,}\\
0 &\text{otherwise.}\end{cases}\end{equation}

\begin{lem}\label{lem:patabran}
Let $K/\mathbb{Q}$ be a non-Galois extension of $\mathbb{Q}$ of degree $3$.
Let $L/\mathbb{Q}$ be the normal closure of $K/\mathbb{Q}$.
Let $p$ be a rational prime that does not ramify in $L/\mathbb{Q}$.
Then $p$ splits completely in $K/\mathbb{Q}$ if and only if it splits
completely in $L/\mathbb{Q}$.
\end{lem}
\begin{proof}
If $\Frob_p = \{I\}$, then $p$ splits completely in $K/\mathbb{Q}$
and in $L/\mathbb{Q}$. Suppose $\Frob_p \ne \{I\}$. Then $p$ does
not split completely in $L/\mathbb{Q}$. Since only one map in
$\Gal(L/\mathbb{Q})$ other than the identity fixes $K/\mathbb{Q}$,
and every conjugacy class in $\Gal(L/\mathbb{Q})$ has more than one
element, there must be a map $\phi \in \Frob_p$ that does not fix
$K/\mathbb{Q}$. Hence $p$ does not split completely in $K/\mathbb{Q}$.
\end{proof}

\begin{lem}\label{lem:taube}
 Let $K/\mathbb{Q}$ be an extension of $\mathbb{Q}$ of degree $3$.
Let $\alpha$ be a positive real number.
Let $S_{\alpha}(X) = \sum_{n\leq X} R(\alpha,n)$.
Then
\begin{equation}\label{eq:atya}\begin{aligned}
 S_{\alpha}(X) &\sim C_{K,\alpha} 
 X (\log X)^{\frac{1}{3} 2^{2 \alpha} - 1} 
     \text{ if $K/\mathbb{Q}$ is Galois,}\\
 S_{\alpha}(X) &\sim C_{K,\alpha} X 
(\log X)^{\frac{1}{2} 2^{\alpha} + \frac{1}{6} 2^{2 \alpha} - 1}
 \text{ if $K/\mathbb{Q}$ is not Galois,}
\end{aligned}\end{equation}
where $C_{K,\alpha}>0$ depends only on $K$ and $\alpha$, and
the dependence of $C_{K,\alpha}$ on $\alpha$ is continuous.
\end{lem}
\begin{proof}
Suppose $K/\mathbb{Q}$ is Galois. Then, for $\Re s>1$,
\[\begin{aligned}
\zeta_{K/\mathbb{Q}}(s) &= \prod_{\mathfrak{p}\in I_K} 
 \frac{1}{1-(N \mathfrak{p})^{-s}}\\ &= \prod_{\text{$p$ ramified}}
 \frac{1}{1- p^{-s}} \mathop{\prod_{\text{$p$ unsplit}}}_{
\text{\& unram.}} \frac{1}{1 - p^{-3 s}}
\prod_{\text{$p$ split}} \frac{1}{(1- p^{-s})^3} .\end{aligned}\]
Hence
\begin{equation}\label{eq:comc}
 \prod_{\text{$p$ split}} (1 + \beta p^{-s}) = L_1(s) 
(\zeta_{K/\mathbb{Q}}(s))^{\beta/3} ,\end{equation}
where $L_1(s)$ is holomorphic and bounded on $\{s : \Re s > 1/2+\epsilon\}$.
Since
\[\sum_n R(\alpha, n) n^{-s} = \mathop{\sum_n}_{p|n \Rightarrow
\text{$p$ split}} 2^{2\alpha} n^{-s} 
= \prod_{\text{$p$ split}} (1 + 2^{2 \alpha} p^{-s} + 2^{2 \alpha} p^{-2 s} +
\dotsb) ,\]
 it follows that
\[\sum_n R(\alpha,n) = L_1(s)
 (\zeta_{K/\mathbb{Q}}(s))^{2^{2 \alpha}/3} .\]
By a Tauberian theorem (see, e.g., \cite{PT}, Main Th.) we can conclude that
\[\frac{1}{X} \sum_{n\leq X} R(\alpha,n)
\sim C_{K,\alpha} (\log X)^{\frac{1}{3} 2^{2 \alpha} - 1}\]
for some positive constant $C_{K,\alpha}>0$.

Now suppose that $K/\mathbb{Q}$ is not Galois. Denote the splitting
type of a prime $p$ in $K/\mathbb{Q}$ by $p = \mathfrak{p}_1 \mathfrak{p}_2$,
$p = \mathfrak{p}_1 \mathfrak{p}_2 \mathfrak{p}_3$, 
$p = \mathfrak{p}_1^2 \mathfrak{p}_2$, etc. Let $L/\mathbb{Q}$ be the
Galois closure of $K/\mathbb{Q}$. By Lemma \ref{lem:patabran},
\begin{equation}\begin{aligned}
\zeta_{K/\mathbb{Q}}(s) &= \prod_{\mathfrak{p}\in I_K} \frac{1}{1 - (N \mathfrak{p})^{-s}} = L_2(s) \prod_{p = \mathfrak{p}_1 \mathfrak{p}_2} 
\frac{1}{(1-p^{-s})} 
\prod_{p = \mathfrak{p}_1 \mathfrak{p}_2 \mathfrak{p}_3}
\frac{1}{(1 - p^{-s})^{3}},\\
\zeta_{L/\mathbb{Q}}(s) &= \prod_{\mathfrak{p}\in I_L} 
 \frac{1}{1 - (N \mathfrak{p})^{-s}} = L_3(s) \prod_{p = \mathfrak{p}_1
\mathfrak{p}_2 \mathfrak{p}_3} \frac{1}{(1 - p^{-s})^6} ,
\end{aligned}\end{equation}
where $L_2(s)$ and $L_3(s)$ are continuous, non-zero
 and bounded on $\{s: \Re s > \frac{1}{2}\}$.
Thus
\[\begin{aligned}
\sum_n R(\alpha,n) n^{-s} &=
 \prod_{p = \mathfrak{p}_1 \mathfrak{p}_2} ( 1 + 2^{\alpha}
p^{-s}) \prod_{p = \mathfrak{p}_1 \mathfrak{p}_2 \mathfrak{p}_3}
 (1 + 2^{2 \alpha} p^{-s}) \\
&= L_4(s) \zeta_{K/\mathbb{Q}}(s)^{2^{\alpha}}
 \zeta_{L/\mathbb{Q}}^{- \frac{1}{2} 2^{\alpha} + \frac{1}{6}
2^{2 \alpha}}. \end{aligned}\]
Since $\zeta_{K/\mathbb{Q}}$ and $\zeta_{L/\mathbb{Q}}$ both
have a pole of order $1$ at $s=1$, we can apply a Tauberian theorem
as before, obtaining
\[\frac{1}{X} \sum_{n\leq X} R(\alpha,n) 
\sim C_{K,\alpha} (\log X)^{\frac{1}{2} 2^{\alpha} +
\frac{1}{6} 2^{2 \alpha} - 1} .\]
\end{proof}
\subsection{The square-free sieve for homogeneous quartics}\label{subs:ludm}
We need the following simple lemma. See Lemma \ref{lem:ramsay} for
a related statement.
\begin{lem}\label{lem:facil}
 Let $f\in \mathbb{Z}\lbrack x,z\rbrack$ be a homogeneous
polynomial. Then there is a constant $C_f$ such that the following holds.
Let $N$ be a positive integer larger than $C_f$. Let $p$ be a 
prime larger than $N$. Then there are at most $12 \deg(f)$ pairs $(x,y)\in 
\mathbb{Z}^2$, $|x|, |z|\leq N$, $\gcd(x,z) = 1$, such that
\begin{equation}\label{eq:maud} p^2 | f(x,z) .\end{equation}
\end{lem}
\begin{proof}
If $N$ is large enough, then $p$ does not divide the discriminant of $f$.
Hence
\begin{equation}\label{eq:bablon}
f(r,1) \equiv 0 \mo p^2\end{equation}
has at most $\deg(f)$ solutions in $\mathbb{Z}/p^2$. If $N$ is large enough for
$p^2$ not to divide the leading coefficients of $f$, then $(x,z) = (1,0)$
does not satisfy (\ref{eq:maud}). Therefore, any solution (x,z) to 
(\ref{eq:maud}) gives us a solution $r = x/z$ to (\ref{eq:bablon}). 
We can focus on solutions $(x,y)\in \mathbb{Z}^2$ to (\ref{eq:maud})
with $x$, $y$ non-negative, as we need only flip signs to repeat the 
procedure for the other quadrants.

Suppose we have two solutions $(x_0,z_0), (x_1, z_1)\in \mathbb{Z}^2$
to (\ref{eq:maud}) such that
\[0\leq |x_0|, |x_1|, |z_0|, |z_1|\leq N,\] 
\[\gcd(x_0,z_0) = \gcd(x_1,z_1) =1 ,\]
\[x_0/z_0 \equiv r \equiv x_1/z_1 \mo p^2.\] Then
\[x_0 z_1 - x_1 z_0 \equiv 0 \mo p^2 .\]
Since $0\leq x_j, z_j\leq N$ and $p>N$, we have that 
\[-p^2 < x_0 z_1 - x_1 z_0 < p^2,\]
and thus $x_0 z_1 - x_1 z_0$ must be zero. Hence $x_0/z_0 = x_1/z_1$.
Since $\gcd(x_0,z_0) = \gcd(x_1,z_1) = 1$ and $\sgn(x_0) = \sgn(x_1)$,
it follows that $(x_0,z_0) = (x_1,z_1)$.
\end{proof} 
\begin{Rem} It was pointed out by Ramsay \cite{Ra} that
an idea akin to that in Lemma \ref{lem:facil} suffices to improve
Greaves's bound for homogeneous sextics \cite{Gr} from
$\delta(N) = N^2 (\log N)^{-1/3}$ to $\delta(N) = N^2 (\log N)^{-1/2}$.
\end{Rem}
\begin{prop}\label{prop:fourcoffees}
Let $f\in \mathbb{Z}\lbrack x,z\rbrack$ be a homogeneous irreducible
polynomial of degree $4$. Let
\[\delta(N) = \{x,z\in \mathbb{Z}^2 : |x|, |z|\leq N, \gcd(x,z)=1,
 \exists p>N \text{ s.t. } p^2 | f(x,y)\} .\]
Then
\[\delta(N)\ll N^{4/3} (\log N)^A,\]
where A and the implied constant depend only on $f$.
\end{prop}
\begin{proof}
Write $A = \max_{|x|, |z|\leq N} f(x,z)$. Clearly $A\ll N^4$. We can
write 
\[\begin{aligned}\delta(N) &\leq \sum_{0 < |d|\leq M} 
\#\{x,y,z\in \mathbb{Z}^3, |x|, |z|\leq N, \gcd(x,z)=1 :
d y^2 = f(x,z) \} \\
&+\sum_{N<p\leq \sqrt{A/M}} \# \{x,z\in \mathbb{Z}^2, |x|, |z|\leq N,
\gcd(x,z) = 1: p^2 | f(x,z) \} .\end{aligned}\]
Let $M\leq N^3$. By Lemma \ref{lem:facil},
\[\sum_{N<p\leq \sqrt{A/M}} \# \{x,z\in \mathbb{Z}^2, |x|, |z|\leq N,
\gcd(x,z) = 1: p^2 | f(x,z) \}\]
is at most a constant times $\frac{1}{\log N} \sqrt{N^{4 - \beta}}$,
where $\beta = (\log M)/(\log N)$. It remains to estimate
\[\sum_{0<|d|\leq M} S(d),\]
where we write
\[S(d) = \#\{x,y,z\in \mathbb{Z}^3, |x|, |z|\leq N, \gcd(x,z)=1 :
d y^2 = f(x,z) \}. \]

Let $C_{f,1}$, $C_{f,2}$, $C_{f,3}$, $C_{f,4}$ be as in Proposition 
\ref{prop:clinton}. Let $K$, $C_f$, $\omega$ and $\omega_K$ be as in
Lemma \ref{lem:notse}. Write $C_{f,5}$ for $C_f$.

By Proposition \ref{prop:clinton},
\[
\sum_{0<|d|< C_{f,4}} S(d) \ll \left(1 + 2 C_{f,3} \sqrt{\frac{7}{2} \log N +
 C_{f,1}}\right)^{C_1} \ll (\log N)^{C_2},\]
where $C_1 = \max_{0<d<C_{f,4}} \rnk(E_d)$, $C_2$ and the implied
constant depend only on $f$. Let $\epsilon$ be a small positive real
number. By Proposition \ref{prop:clinton} and Lemma \ref{lem:notse},
\[\begin{aligned}
\sum_{C_{f,4}\leq |d|<N^{\epsilon}} S(d) &\ll
\sum_{C_{f,4}\leq |d|<N^{\epsilon}} \left(1 + 2 \sqrt{\frac{7}{2} \log N +
C_{f,1}}\right)^{\rnk(E_d)} \\
&\ll \sum_{C_{f,4}\leq |d|<N^{\epsilon}} \left(1 + 2 \sqrt{\frac{7}{2} \log N +
C_{f,1}}\right)^{C_{f,5} + \omega_K(d) - \omega(d)}.\end{aligned}\]
We have the following crude bounds:
\begin{equation}\label{eq:cru}
\omega(d) \leq \frac{\log |d|}{\log \log |d|},\;\;
\omega_K(d)\leq 3 \omega(d) .\end{equation}
Hence
\[\begin{aligned}
\sum_{C_{f,4}\leq |d|<N^{\epsilon}} 
S(d) &\ll
\sum_{C_{f,4}\leq d<N^{\epsilon}} (\log N)^{C_{f,5} + 2 \log d/\log \log d}
\\ &\leq
N^{\epsilon} (\log N)^{C_1} (\log N)^{2 \epsilon \log N/\log \log N} 
\leq (\log N)^{C_1} N^{3 \epsilon} ,\end{aligned}\]
where $C$ depends only on $f$ and $\epsilon$.
For any $d$ with $|d|>N^{\epsilon}$, 
 Proposition \ref{prop:clinton} and Lemma \ref{lem:notse} give us
\[\begin{aligned}
S(d) 
&\ll \left(1 + 2 \sqrt{\left(\frac{7}{2} \log N + C_{f,1}\right)/\left(\frac{1}{8} 
 \epsilon \log N + C_{f,2}\right)}\right)^{\rnk(E_d)} \\
&\ll (12 \epsilon^{-1/2})^{C_{f,5} + \omega_K(d) - \omega(d)} 
\leq 2^{C_2 \omega_K(d) - C_2 \omega_K(d)} ,\end{aligned}\]
where $C_2$ depends only on $f$ and $\epsilon$. By Lemma \ref{lem:taube}
we can conclude that 
\[
\sum_{N^{\epsilon}<|d|\leq M} S(d) \ll \sum_{d=1}^M 2^{C_2
\omega_K(d) - C_2 \omega_K(d)} \ll C_3 M (\log N)^{C_4},\]
where $C_3$ and $C_4$ depend only on $f$ and $\epsilon$. Set $M = N^{4/3}$,
$\epsilon = 1/4$.
\end{proof}

\subsection{Homogeneous cubics}

\begin{prop}\label{prop:aulait}
Let $f\in \mathbb{Z}\lbrack x,z\rbrack$ be a homogeneous irreducible
polynomial of degree $3$. Let
\[\delta(N) = \{x,z\in \mathbb{Z}^2 : |x|, |z|\leq N, \gcd(x,z)=1,
 \exists p>N \text{ s.t. } p^2 | f(x,y)\} .\]
Then
\[\delta(N)\ll N^{4/3} (\log N)^A,\]
where A and the implied constant depend only on $f$.
\end{prop}
\begin{proof}
Write $A = \max_{|x|, |z|\leq N} f(x,z)$. Clearly $A\ll N^4$. We can
write 
\[\begin{aligned}\delta(N) &\leq 
\sum_{0 < |d|\leq M} 
\#\{x,y,z\in \mathbb{Z}^3, |x|, |z|\leq N, \gcd(x,z)=1 :
d y^2 = f(x,z) \} \\
&+\sum_{N<p\leq \sqrt{A/M}} \# \{x,z\in \mathbb{Z}^2, |x|, |z|\leq N,
\gcd(x,z) = 1: p^2 | f(x,z) \} .\end{aligned}\]
Let $M\leq N^2$. 
By Lemma \ref{lem:facil}, the second term on the right is at most
a constant times
$N^{2-\beta/2}/\log N$. Now notice that any point
$(x,y,z)\in \mathbb{Z}^3$ on $d y^2 = f(x,z)$ gives us a rational
point $(x',y') = (x/z,y/z^2)$ on
\begin{equation}\label{eq:tineq}
d' {y'}^2 = f(x',1),\end{equation}
where $d' = d z$. Moreover, a rational point on (\ref{eq:tineq})
can arise from at most one point $(x,y,z)\in \mathbb{Z}^3$, $\gcd(x,z)=1$,
in the given fashion. 

If $d\leq M$, then $|d'|=|d z|\leq M N$. The height $h_x(P)$
of the point $P=(x/z,y/z^2)$ is at most $N$. It follows by Lemma 
\ref{lem:soide} that $\hat{h}(P)\leq N + C_f$, where $C_f$ is a
constant depending only on $f$. By Corollaries \ref{cor:nosilv} and
\ref{cor:magd}, there are at most
\[O\left((1 + 2 \sqrt{(\log N + C_f')/(\log |d| + C_f)})^{\rnk(E_d)}\right)\]
rational points $P$ of height $\hat{h}(P) \leq N + C_f$.
We proceed as in Proposition \ref{prop:fourcoffees}, and obtain that
\[\sum_{0 < |d|\leq M} 
\#\{x,y,z\in \mathbb{Z}^3, |x|, |z|\leq N, \gcd(x,z)=1 :
d y^2 = f(x,z) \}\]
is at most $O(M N (\log N))^A$. Set $\beta = 1/3$.
\end{proof}
\subsection{Homogeneous quintics}
We extract the following result from \cite{Gr}. 
\begin{lem}\label{lem:tumba}
Let $f\in \mathbb{Z}\lbrack x,y\rbrack$ be a homogeneous irreducible
polynomial of degree at most $5$. For all $M<N^{\deg f}$, $\epsilon>0$,
\[
\sum_{d=1}^M \#\{x,y,z\in \mathbb{Z}^3, |x|, |z|\leq N, \gcd(x,z)=1 :
d y^2 = f(x,z) \}\]
is at most a constant 
times $N^{(18-\frac{1}{2} \beta^2)/(10-\beta) + \epsilon}$,
where $\beta = (\log M)/(\log N)$. The implied constant depends only on $f$ 
and $\epsilon$.
\end{lem}
\begin{proof}
By \cite{Gr}, Lemmas 5 and 6, where the parameters $d$ and $z$ (in
the notation of \cite{Gr}) are set to the values $d=1$ and
 $z = N^{(1-\beta/2)/(5/2 - \beta/4)}$.
\end{proof}
\begin{prop}\label{prop:jocha}
Let $f\in \mathbb{Z}\lbrack x,z\rbrack$ be a homogeneous irreducible 
polynomial of degree $5$. Let
\[\delta(N) = \{x,z\in \mathbb{Z}^2 : |x|, |z|\leq N, \gcd(x,z)=1,
 \exists p>N \text{ s.t. } p^2 | P(x,y)\} .\]
Then, for any $\epsilon>0$,
\[\delta(N) \ll 
N^{(5+\sqrt{113})/8 + \epsilon}\]
where the implied constant depends only on $f$ and $\epsilon$.
\end{prop}
\begin{proof}
Let $A = \max_{|x|, |z|\leq N} f(x,z)$. Clearly $A\ll N^{\deg(f)}$.
We can write
\[\begin{aligned}
\delta(N) &\leq \sum_{0<|d|\leq M} 
\#\{x,y,z\in \mathbb{Z}^3, |x|, |z|\leq N, \gcd(x,z)=1 :
d y^2 = f(x,z) \} \\
&+\sum_{N<p\leq \sqrt{A/M}} \# \{x,z\in \mathbb{Z}^2, |x|, |z|\leq N,
\gcd(x,z) = 1: p^2 | f(x,z) \} .\end{aligned}\]
By Lemmas \ref{lem:tumba} and \ref{lem:facil}, 
\[\delta(N) \ll N^{(18 - \frac{1}{2} \beta^2)/(10 - \beta) + \epsilon} +
\frac{1}{\log N} \sqrt{N^{\deg(f) - \beta}},\]
where $\beta = (\log M)/(\log N)$. Set $\beta = (15 - \sqrt{113})/4$.
\end{proof}

\subsection{Quasiorthogonality, kissing numbers and cubics}\label{subs:kiss}
\begin{lem}\label{lem:agamen}
Let $f\in \mathbb{Z}\lbrack x\rbrack$ be a cubic polynomial of non-zero
discriminant. Let $d$ be a square-free integer. Then, for any two 
distinct
integer points $P=(x,y)\in \mathbb{Z}^2$, $P'=(x',y')\in \mathbb{Z}^2$
on the elliptic curve
\[E_d : d y^2 = f(x),\]
we have
\[\hat{h}(P+P') \leq 3 \max(\hat{h}(P),\hat{h}(P')) + C_f,\]
where $C_f$ is a constant depending only on $f$.
\end{lem}
Cf. \cite{GS}, Proposition 5.
\begin{proof}
Write $f(x) = a_3 x^3 + a_2 x^2 + a_1 x + a_0$. Let $P + P' = (x'',y'')$.
By the group law,
\[\begin{aligned}
x'' &= \frac{d (y_2 - y_1)^2}{a_3 (x_2 - x_1)^2} - \frac{a_2}{a_3} - x_1 - x_2
\\
&= \frac{d (y_2 - y_1)^2 - a_2 (x_2 - x_1)^2 - a_3 (x_2 - x_1)^2 (x_1 + x_2)}{
a_3 (x_2 - x_1)^2} .\end{aligned}\]
Clearly $|a_3 (x_2 - x_1)^2|\leq 4 |a_3| \max(|x_1|^2,|x_2|^2)$. Now
\[|d (y_2 - y_1)^2| \leq 4 |d| \max(y_1^2, y_2^2) = 4 \max(|f(x_1)|,|f(x_2)|) .\]
Hence
\[|d (y_2 - y_1)^2 - a_2 (x_2 - x_1)^2 - a_3 (x_2 - x_1)^2 (x_1 + x_2)|
 \leq A \max( |x|^3, |x'|^3) ,\]
where $A$ is a constant depending only on $f$. Therefore
\[\begin{aligned}
h_x(P) &= \log(\max(|\num(x'')|,|\den(x'')|)) \\
 &\leq 3 \max(\log |x|, \log |x'|) + \log A \\
 &\leq 3 \max(h_x(P),h_x(P')) + \log A .\end{aligned}\]
By Lemma \ref{lem:soide}, the difference $|\hat{h}-h_x|$ 
is bounded by a constant
independent of $d$. The statement follows immediately.
\end{proof}
Consider the elliptic curve
\[E_d:d y^2 = f(x) .\]
There is a $\mathbb{Z}$-linear map from $E_d(\mathbb{Q})$ to
$\mathbb{R}^{\rnk(E_d)}$ taking the square of
the Euclidean norm back to the canonical height. In other words, any
 given integer point $P=(x,y)\in E_d$ will be taken
to a point $L(P)\in \mathbb{R}^{\rnk(E_d)}$ whose Euclidean norm $|L(P)|$
satisfies
\[|L(P)|^2 = \hat{h}(P) = \log x + O(1) ,\]
where the implied constant depends only on $f$. In particular, the set of all
integer points $P=(x,y)\in E_d$ with
\begin{equation}\label{eq:perdie}
N^{1-\epsilon}\leq x\leq N\end{equation}
will be taken to a set of points $L(P)$ in $\mathbb{R}^{\rnk(E_d)}$ with
\[(1-\epsilon) \log N + O(1) \leq |L(P)|^2 \leq \log N + O(1) .\]
Let $P, P'\in E_d$ be integer points satisfying (\ref{eq:perdie}).
Assume $L(P)\ne L(P')$. By Lemma \ref{lem:agamen},
\[|L(P) + L(P')|^2 = |L(P + P')|^2 \leq 3 \max(|L(P)|^2,|L(P')|^2) + O(1).\]
Therefore, the inner product $L(P) \cdot L(P')$ satisfies
\[\begin{aligned} L(P) \cdot L(P') &= \frac{1}{2} (|L(P) + L(P')|^2 -
(|L(P)|^2 + |L(P')|^2)) \\
&\leq \frac{1}{2} (3 \max(|L(P)|^2,|L(P')|^2) + O(1) - (|L(P)|^2 +
|L(P')|^2))\\
&\leq \frac{1}{2} ((1+\epsilon) \log(N) + O(1)) \\
&\leq \frac{1}{2} \frac{(1+\epsilon) + O((\log N)^{-1})}{(1-\epsilon)^2}
 |L(P)| |L(P')| .\end{aligned}\]
We have proven
\begin{lem} 
Let $f\in \mathbb{Z}\lbrack x\rbrack$ be a cubic polynomial
of non-zero discriminant. Let $d$ be a square-free integer. Consider the
elliptic curve
\[E_d : d y^2 = f(x) .\]
Let $S$ be the set
\[\{(x,y)\in \mathbb{Z}^2 : N^{1-\epsilon} \leq |x|\leq N, d y^2 = f(x)\} .\]
Let $L$ be a linear map taking $E(\mathbb{Q})$ to $\mathbb{R}^{\rnk(E_d)}$
and the square of the Euclidean norm back to
the canonical height $\hat{h}$.
Then, for any distinct points $P, P' \in L(S) \subset \mathbb{R}^{\rnk(E_d)}$
with the angle $\theta$ between $P$ and $P'$ is at least
\[\arccos\left(\frac{1}{2} \frac{(1 + \epsilon) + O((\log N)^{-1})}{(1-\epsilon)^2} \right) 
= 60^\circ + O(\epsilon + (\log N)^{-1}) ,\]
where the implied constant depends only on $f$.
\end{lem}

Let $A(\theta,n)$ be the maximal number of points that
can be arranged in $\mathbb{R}^n$ with angular separation no smaller
than $\theta$. Kabatiansky and Levenshtein (\cite{KL}; vd. also
\cite{CS}, (9.6)) show that, for $n$ large enough,
\begin{equation}\label{eq:kale}\frac{1}{n} \log_2 A(n,\theta) \leq
\frac{1 + \sin \theta}{2 \sin \theta} \log_2 
 \frac{1 + \sin \theta}{2 \sin \theta} -
\frac{1 - \sin \theta}{2 \sin \theta} \log_2 
 \frac{1 - \sin \theta}{2 \sin \theta} .\end{equation}
Thus we obtain
\begin{cor}\label{cor:coeur}
Let $f\in \mathbb{Z}\lbrack x\rbrack$ be a cubic polynomial
of non-zero discriminant. Let $d$ be a square-free integer. Consider the
elliptic curve
\[E_d : d y^2 = f(x) .\]
Let $S$ be the set
\[\{(x,y)\in \mathbb{Z}^2 : N^{1-\epsilon} \leq |x|\leq N, d y^2 = f(x)\} .\]
Then
\[\# S \ll 2^{(\alpha + O(\epsilon + (\log N)^{-1})) \rnk(E_d)},\]
where 
\begin{equation}\label{eq:calor}
\alpha = \frac{2 + \sqrt{3}}{2 \sqrt{3}} \log_2
\frac{2 + \sqrt{3}}{2 \sqrt{3}} - \frac{2 - \sqrt{3}}{2 \sqrt{3}} \log_2
\frac{2 - \sqrt{3}}{2 \sqrt{3}} = 0.4014\dotsc\end{equation}
and the implied constants depend only on $f$.
\end{cor}
Notice that we are using the fact that the size of the torsion group
is bounded.
\begin{prop}\label{prop:bethoo}
Let $f\in \mathbb{Z}\lbrack x\rbrack$ be an irreducible cubic polynomial.
Let
\[\delta(N) = \{1\leq x \leq N : \exists p > N^{1/2} \text{ s.t. }
p^2 | f(x) \} .\]
Then
\begin{equation}\label{eq:gerg}\delta(N) \ll N (\log N)^{-\beta},
\end{equation}
where
\[\beta = 1 - \frac{1}{9} 2^{2 \alpha} = 0.8061\dotsc\] 
if the discriminant of $f$ is a square, 
\[\beta = 1 - \frac{1}{6} 2^{\alpha} - \frac{1}{18} 2^{2 \alpha} =
0.6829\dotsc\]
if the discriminant of $f$ is not a square, and $\alpha$ is as in
(\ref{eq:calor}).
The implied constant in (\ref{eq:gerg}) depends only on $f$.
\end{prop}
\begin{proof}
Let $A = \max_{1\leq x\leq N} f(x)$. Clearly $A\ll N^3$. We can write
\[\begin{aligned}
\delta(N) &\leq \sum_{N^{1/2}<p<\sqrt{A/M}} \# \{1\leq x\leq N : p^2 | f(x)\}\\
&+ \# 
\{1\leq x\leq N^{1-\epsilon} : \exists p>N^{1/2} \text{ s.t. }
p^2 | f(x) \} \\
&+ \sum_{1\leq |d|\leq M} \#\{x,y\in \mathbb{Z}^2 : N^{1-\epsilon}\leq x\leq N,
d y^2 = f(x) \} .\end{aligned}\]
Let $M\leq N^2$. Then
the first term is at most
\[\sum_{N^{1/2} < p < \sqrt{A/M}} 3\ll \frac{3 \sqrt{A/M}}{\log \sqrt{A/M}} \ll \frac{N^{3/2} M^{-1/2}}{\log N} .\]
The second term is clearly no greater than $N^{1-\epsilon}$.
It remains to bound
\[\sum_{1\leq |d|\leq M} B(d),\]
where \[B(d) = 
\#\{x,y\in \mathbb{Z}^2 : N^{1-\epsilon}\leq x\leq N,
d y^2 = f(x) \} .\]
By Lemma \ref{lem:notse}, Corollary \ref{cor:coeur}
and the remark before (\ref{eq:manandur})
\[B(d) \ll R(\alpha+ O(\epsilon + (\log N)^{-1}),d),\]
where $K$ is as in Lemma \ref{lem:notse}
and $\alpha$ is as in Corollary \ref{cor:coeur} and $R(\alpha,d)$ is
as in (\ref{eq:manandur}). Thanks to (\ref{eq:cru}),
we can omit the term $O((\log N)^{-1})$ from the exponent.
Hence it remains to estimate
$S(M) =
\sum_{1\leq d\leq M} R(\alpha + O(\epsilon + (\log N)^{-1}),d)$.
By Lemma \ref{lem:taube},
\[\begin{aligned}
 S(M) &\ll
 M (\log M)^{\frac{1}{3} 2^{2(\alpha+O(\epsilon))}-1} 
      \text{ if $K/\mathbb{Q}$ is Galois,}\\
 S(M) &\ll
 M (\log M)^{\frac{1}{2} 2^{\alpha+O(\epsilon)} +
\frac{1}{6} 2^{2 (\alpha+O(\epsilon))} - 1} 
 \text{ if $K/\mathbb{Q}$ is not Galois.}\end{aligned}\]
Set
\[\begin{aligned}
M &= N (\log N)^{- \frac{2}{9} 2^{2\alpha}} \text{ if $K/\mathbb{Q}$ is Galois}\\
M &= N (\log N)^{- \frac{1}{3} 2^{\alpha} - \frac{1}{9} 2^{2 \alpha}}
 \text{ if $K/\mathbb{Q}$ is not Galois.}\end{aligned}\]
Let $\epsilon = (\log \log M)^{-1}$. Since $K/\mathbb{Q}$ is
Galois if and only if the discriminant of $f$ is a square, the statement
follows.
\end{proof}
\subsection{Hyperelliptic curves, Mumford's gap and sextics}
We will first need a few lemmas on hyperelliptic curves.
 See \cite{BK}, pp. 718--734, for the analogous
statements on elliptic curves. Cf. also \cite{CF}, \S 7.5, and \cite{FPS},
Prop. 3. As is usual, when we speak of the curve $C:y^2 = f(x)$, 
$\deg(f)\geq 5$, we mean the non-singular projective curve in which the affine
curve $y^2 = f(x)$ is contained. 

\begin{lem}\label{lem:aghar}
Let $f(x) \in K\lbrack x\rbrack$ be a polynomial of degree $6$ defined over
a field $K$ of characteristic zero. Let $C$ be the curve $y^2 = f(x)$.
 Let $A_K$ be the $K$-algebra
$K\lbrack x\rbrack/f\lbrack x\rbrack$. Write $A_K = \oplus_j K_j$, where
each $K_j$ is a finite extension of $K$ corresponding to an irreducible
factor of $K$. Define $N_{A_K/K}:A_K\to K$ to be the product of all norms
$N_{K_j/K}:K_j\to K$. Let $J$ be the Jacobian of $C$. Then there
is a homomorphism
\[\mu:J(K) \to A_K^*/K^* (A_K^*)^2\]
whose image is in the kernel of
\[A^*/K^*(A_K^*)^2 \xrightarrow[N_{A_K/K}]{} K^*/(K^*)^2\]
and whose kernel is generated by $2 J(K)$ and at most one element of
$J(K)$ not in $2 J(K)$. The map $\mu$ commutes with field extensions
$K\to K'$:
\begin{equation}\label{eq:commie}\begin{CD}
J(K) @> >> A_K^*/K^* (A_K^*)^2 \\
@ VV V @ VV V \\
J(K') @> >> A_{K'}^*/K'^* (A_{K'}^*)^2 \end{CD}\end{equation}
\end{lem}
\begin{proof}
See \cite{Ca}, p. 30, for the statement on the image and kernel of $\mu$.
The commutativity of (\ref{eq:commie}) follows from the construction of
$\mu$ in \cite{Ca}, pp. 50--51.
\end{proof}
Given $K$ and $A_K = \oplus_j K_j$ as in Lemma \ref{lem:aghar}, we speak
of the coordinates $(x_j) \in \oplus_j K_j$ of a point $x\in A_K$ in the
natural fashion. 
\begin{lem}\label{lem:zolt}
Let $f(x) \in K\lbrack x\rbrack$ be a polynomial of degree $6$ defined over
a $\mathfrak{p}$-adic field $K$ whose residue field has odd characteristic.
Suppose $C: y^2 = f(x)$ has good reduction. Let $\mu$ 
be as in Lemma \ref{lem:aghar}. Let $(x_j)$ be the
coordinates of a representative $x\in A_K^*$ of a point $\mu(y) \in
A_K^*/K^* (A_K^*)^2$ in the image of $\mu$. Then $v_{K_j}(x_j) \mod 2$
is independent of $j$.
\end{lem}
\begin{proof}
Let $k$ be the residue field of $K$, $\tilde{K}$ the maximal unramified
extension of $K$, and $\bar{k}$ the algebraic closure of $k$. Let 
$J_1(\tilde{K})$ be the kernel of the residue map $J(\tilde{K}) \to
J(\bar{k})$. We have an exact sequence
\[0\to J_1(\tilde{K}) \to J(\tilde{K}) \to J(\bar{k}) \to 0,\]
from which we obtain
 \[J_1(\tilde{K})/2 J_1(\tilde{K}) \to 
J(\tilde{K})/2 J(\tilde{K}) \to J(\bar{k})/2 J(\bar{k}) \to 0.\]
Since $\bar{k}$ is algebraically closed, the group $J(\bar{k})/2 J(\bar{k})$
is trivial. Since $J_1(\tilde{K})$ is isomorphic to a formal group over $k$
(\cite{CF}, (7.3.5)) and the characteristic of $k$ is not $2$, the group 
$J(\tilde{K})/2 J(\tilde{K})$ is trivial as well.  Hence 
$J(\tilde{K})/2 J(\tilde{K})$ is trivial.

By virtue of the commutative diagram
\[ \begin{CD} 
       J(K)/2 J(K) @>{\mu}>> A_K^*/K^* (A_K^*)^2\\
@VV V @VV V \\
J(\tilde{K})/2 J(\tilde{K}) = \{e\} @>{\mu}>> 
 A_{\tilde{K}}^*/\tilde{K}^* (A_{\tilde{K}}^*)^2,\\
       \end{CD} 
       \] 
the image of $\mu(J(K)/2 J(K))$ lies in $A_{\tilde{K}}^* \cap 
\tilde{K}^* (A_{\tilde{K}}^*)^2$. 
Since, for every index $j$, the extensions $\tilde{K}_j/
K_j$ and $\tilde{K}_j/\tilde{K}$ are unramified,
we can write, for any $r\in \tilde{K}^*$,
$s\in A_{\tilde{K}}^*$ with $r s^2 \in A_{\tilde{K}}^*$,
\[ v_{K_j}(r s^2) = v_{\tilde{K}_j}(r s^2) =  v_{\tilde{K}_j}(r)
+ 2 v_{\tilde{K}_j}(s) = v_{\tilde{K}}(r) + 2 v_{\tilde{K}_j}(s) .\]
The statement follows.
\end{proof}

\begin{prop}\label{prop:malik}
Let $f(x) \in \mathbb{Q}\lbrack x\rbrack$ be an irreducible polynomial of
degree $6$. Choose a root $x_0$ of $f(x)=0$. Let $L = \mathbb{Q}(x_0)$.
For every square-free integer $d$, let $C_d$ be the hyperelliptic curve
$d y^2 = f(x)$. Then
\[\rnk(C_d) \leq \kappa_f + \sum_{p|d} (n_p - m_p),\]
where \[\begin{aligned}
n_p &= \text{number of primes of $L$ above $p$,}\\
m_p &= \begin{cases} 1 & \text{if all primes above $p$ are of even inertia
degree,} \\ 2 & otherwise, \end{cases} \end{aligned}\]
and $\kappa_f$ is a constant depending only on $f$.
\end{prop}
\begin{proof}
We can identify $L$ with the algebra $A_{\mathbb{Q}}$ from Lemmas
\ref{lem:aghar}--\ref{lem:zolt}; notice that $L$ and $A_{\mathbb{Q}}$ 
are independent of the twist $d$. Consider the exact sequence
\[0\longrightarrow H\longrightarrow A_{\mathbb{Q}}^*/\mathbb{Q}^*
(A_{\mathbb{Q}}^*)^2 \longrightarrow \prod_{\nu} 
(A_{\mathbb{Q}_\nu}^*/\mathbb{Q}_{\nu}^* (A_{\mathbb{Q}_\nu}^*)^2)/
 U_{A_{K_\nu}},\]
where $U_{A_{K_{\nu}}}$ is 
$V_{A_{K_\nu}} \mathbb{Q}_\nu^* (A_{\mathbb{Q}_\nu}^*)^2$,
$V_{A_{K_\nu}}$
is the subgroup of $A_{\mathbb{Q}_\nu}^*$
consisting of the elements $x\in A_{\mathbb{Q}_\nu}^*$ all of whose
coordinates are units, and $H$ is the subgroup of
$L^*/\mathbb{Q}^* (L^*)^2 = A_{\mathbb{Q}}^*/\mathbb{Q}^* (A_{\mathbb{Q}}^*)^2$
consisting of such cosets as are represented by elements
$x\in L^*$ with $(x) = r \mathfrak{s}^2$ for some $r\in \mathbb{Q}$,
$\mathfrak{s}\in I_L$. Clearly
\[\rnk(H) \leq 5 + \rnk(h_L/h_L^2),\]
where $h_L$ is the class group of $L$. Hence the rank of the kernel
of the composition
\[c: J_{C_d}(\mathbb{Q})/2 J_{C_d}(\mathbb{Q}) 
\xrightarrow{\mu} A_{\mathbb{Q}}^*/\mathbb{Q}^* (A_{\mathbb{Q}}^*)^2
 \rightarrow \prod_{\nu} (A_{\mathbb{Q}_\nu}^*/\mathbb{Q}_{\nu}^*
(A_{\mathbb{Q}_\nu}^*)^2)/U_{A_{\mathbb{Q}_{\nu}}}\]
is at most $6+ \rnk(h_L/h_L^2)$. We have a commutative diagram
\[ \begin{CD} 
       J_{C_d}(\mathbb{Q})/2 J_{C_d}(\mathbb{Q}) @>{c}>> 
\prod_{\nu} (A_{\mathbb{Q}_\nu}^*/\mathbb{Q}_{\nu}^*
(A_{\mathbb{Q}_\nu}^*)^2)/U_{A_{\mathbb{Q}_{\nu}}} 
\\
@VV V @VV V \\
J_{C_d}(\mathbb{Q}_{\nu})/2 J_{C_d}(\mathbb{Q}_{\nu}) @>{\mu}>> 
 (A_{\mathbb{Q}_\nu}^*/\mathbb{Q}_{\nu}^*
(A_{\mathbb{Q}_\nu}^*)^2)/U_{A_{\mathbb{Q}_{\nu}}} .
       \end{CD} 
       \] 
By Lemma \ref{lem:zolt}, the bottom row has trivial image when
$\nu \ne 2,\infty$ and $C_d$ has good reduction at $\nu$. Thus
\[\begin{aligned}
\rnk(J_{C_d}(\mathbb{Q})/2 J_{C_d}(\mathbb{Q})) &\leq 6 + \rnk(h_L/h_L^2) \\
&+ 
\mathop{\mathop{\sum_{\nu = \infty, 2 \text{ or}}}_{\text{L ramifies at $\nu$  
or}}}_{\text{$C_1$ has bad red. at $\nu$}} \rnk(
(A_{\mathbb{Q}_\nu}^*/\mathbb{Q}_{\nu}^*
(A_{\mathbb{Q}_\nu}^*)^2)/U_{A_{\mathbb{Q}_{\nu}}} \\ &+
\mathop{\sum_{p|d}}_{\text{$p$ unram. in $L/\mathbb{Q}$}}
\rnk(\mu(J_{C_d}(\mathbb{Q})/2 J_{C_d}(\mathbb{Q}))) .\end{aligned}\]
By Lemma \ref{lem:aghar}, $\mu(J_{C_d}(\mathbb{Q})/2 J_{C_d}(\mathbb{Q}))$
is contained in the intersection of 
the kernel of $N_{A_\mathbb{Q}/\mathbb{Q}}$ and
$(A_{\mathbb{Q}_\nu}^*/\mathbb{Q}_{\nu}^*
(A_{\mathbb{Q}_\nu}^*)^2)/U_{A_{\mathbb{Q}_{\nu}}}$.
Hence the image of \[\mu(J_{C_d}(\mathbb{Q}_\nu)/2 J_{C_d}(\mathbb{Q}_\nu))\]
under the natural injection
\[(A_{\mathbb{Q}_\nu}^*/\mathbb{Q}_{\nu}^*
(A_{\mathbb{Q}_\nu}^*)^2)/U_{A_{\mathbb{Q}_{\nu}}} \longrightarrow 
(\oplus_{\mathfrak{p}|p} \mathbb{Z}/2 \mathbb{Z})/(1,\dotsc,1)\]
consists of cosets $(x_{\mathfrak{p}})_{\mathfrak{p}|p} \in
 \oplus_{\mathfrak{p}|p} \mathbb{Z}/2 \mathbb{Z}$ with 
$\sum_{\mathfrak{p}|p} f_{\mathfrak{p}} x_{\mathfrak{p}}$ even. There
are $2^{n_p - m_p}$ such cosets, where $n_p$ and $m_p$ are as in the
statement of the present lemma. Therefore
\[\begin{aligned}
\rnk(J_{C_d}(\mathbb{Q}))&\leq \rnk(J_{C_d}(\mathbb{Q})/
2 J_{C_d}(\mathbb{Q})) \leq 6 + \rnk(h_L/h_L^2) \\ &+ 
\mathop{\mathop{\sum_{\nu = \infty, 2 \text{ or}}}_{\text{at $\nu$, $L$ ramifies  
or}}}_{\text{$C_1$ has bad red.}} \rnk(
(A_{\mathbb{Q}_\nu}^*/\mathbb{Q}_{\nu}^*
(A_{\mathbb{Q}_\nu}^*)^2)/U_{A_{\mathbb{Q}_{\nu}}} + \sum_{p|d} (n_p - m_p) .\end{aligned}\]
\end{proof}

\subsection{More averages of divisor functions}
The following lemma is an analytical variant on Chebotarev's theorem. 
We denote the conjugacy
class of a map $\phi \in \Gal(K/\mathbb{Q})$  by $\langle\phi\rangle$.  
\begin{lem}\label{lem:chebo}
Let $K/\mathbb{Q}$ be a Galois number field. Let $\langle\gamma\rangle$  
be a conjugacy class
in $G=\Gal(K/\mathbb{Q})$. Then, for $\Re(s)>1$,
\[\mathop{\prod_p}_{\Frob_p = \langle\gamma\rangle} (1 - p^{-s})^{-1} = 
L_0(s) \prod_\chi L_{\chi}(s)^{\overline{\chi(\gamma)} \cdot \# \langle
\gamma \rangle / \# G},\]
where \begin{itemize}
\item $L_0$ is holomorphic and bounded 
on $\{s\in \mathbb{C}: \Re s \geq 1/2 + \epsilon\}$,
\item the product $\prod_{\chi}$ is taken over all characters $\chi$
of $G$, 
\item $L_{\chi}$ is the Artin $L$-function corresponding
to the character $\chi$. 
\end{itemize}
\end{lem}
\begin{proof}
For every character $\chi$ of $G$, the $L$-function $L_{\chi}$ is of the
form \[L_0(s) \cdot \prod_p (1 - \chi(\Frob_p) p^{-s})^{-1},\]
where $L_0(s)$ is some function holomorphic and bounded for $\Re s \geq
1/2 + \epsilon$. By the orthogonality of the characters of $G$, we have that
 $\sum_{\chi} \chi(\Frob_p) \overline{\chi(\gamma)}$
is $|G|$ for $\langle \Frob_p \rangle = \langle \gamma \rangle$, and
$0$ otherwise.
\end{proof}
\begin{cor}\label{cor:evaine}
 Let $K/\mathbb{Q}$ be a Galois number field. Let $g:\Gal(K/\mathbb{Q})
\to \mathbb{C}$ be given.
Then
\[\sum_{d\leq X} \prod_{p|d} g(\Frob_p) =
C_{K,g} X (\log X)^{-1 + \frac{1}{\# G} \sum_{\gamma\in G} g(\gamma)},\]
where $G = \Gal(K/\mathbb{Q})$, $C_{K,g}$ is a positive constant
depending only on $K$ and $g(\gamma)$ for all $\gamma\in G$,
and the dependence of $C_{K,g}$ on $g(\gamma)$ is continuous for
every $\gamma \in G$.
\end{cor}
\begin{proof}
Clearly
\[\sum_{n} \prod_{p|d} g(\Frob_p) n^{-s} =
 \prod_{\langle \gamma \rangle} 
\mathop{\prod_{p}}_{\Frob_p = 
\langle \gamma \rangle}
  (1 + g(\gamma) p^{-s} + g(\gamma) p^{-2 s} + \dotsb),\]
where the product $\prod_{\langle \gamma \rangle}$ is taken
over all conjugacy classes $\langle \gamma \rangle$ in $G$.
By Lemma \ref{lem:chebo},
\[
\mathop{\prod_{p}}_{\Frob_p = 
\langle \gamma \rangle}
  (1 + g(\gamma) p^{-s} + g(\gamma) p^{-2 s} + \dotsb) = 
L_0(s) \prod_\chi L_{\chi}(s)^{g(\gamma)
\overline{\chi(\gamma)} \cdot \# \langle
\gamma \rangle / \# G},\]
where $L_0$ is holomorphic and bounded on $\Re s \geq 1/2 + \epsilon$.
Hence
\[\sum_{n} \prod_{p|d} g(\Frob_p) n^{-s} = L_0(s) 
 \prod_{\chi} L_{\chi}(s)^{\sum_{\langle g\rangle}
g(\gamma) \overline{\chi(\gamma)} \cdot \# \langle
\gamma \rangle / \# G} .\]
Artin $L$-functions associated to non-principal characters $\chi$
are holomorphic and bounded on a neighbourhood of $s=1$. By a
Tauberian theorem (e.g., \cite{PT}, Main Th.), the statement follows.
\end{proof}
Let $K/\mathbb{Q}$ be a number field; let $K'/\mathbb{Q}$ be its Galois closure. We can see $\Gal(K'/\mathbb{Q})$ as a transitive
permutation group on
$\{1,2,\dotsc,\deg(K/\mathbb{Q})\}$.
Since every permutation of $\{1,2,\dotsc,n\}$ is a product of
disjoint cycles,
 we may speak of the cycles of a map $\gamma\in \Gal(K'/\mathbb{Q})$.
We write $G_1(K'/\mathbb{Q})$ for the set of all $\gamma\in \Gal(K'/\mathbb{Q})$ fixing at least one point in $\{1,2,\dotsc,\deg(K,\mathbb{Q})\}$

\begin{cor}\label{cor:fasc}
Let $f(x) \in \mathbb{Q}\lbrack x\rbrack$ be an irreducible polynomial of
degree $6$. For every square-free integer $d$, let $C_d$ be the hyperelliptic
curve $d y^2 = f(x)$. Let $\alpha$ be a positive real number. 
Write $K'$ for the splitting field of $f$.
Choose a root $x_0$ of $f(x)=0$. Let $K = \mathbb{Q}(x_0)$.
Let $R(\alpha,d) = 2^{\alpha \rnk(C_d)}$ if no $p|d$ is unsplit in
$K/\mathbb{Q}$. Let $R(\alpha,d) = 0$ otherwise. 
Then
\begin{equation}\label{eq:aght}
\sum_{d\leq X} R(\alpha,d) \ll X (\log X)^{\delta-1},\end{equation}
where
\[\delta = \frac{1}{\Gal(K'/\mathbb{Q})} \sum_{\gamma \in G_1(K'/\mathbb{Q})}
 2^{\alpha (n_{\gamma} - 2)},\] $n_{\gamma}$ is the number
of cycles of $\gamma$, 
and the implied constant in (\ref{eq:aght}) depends only on $K$ and
$\alpha$; the dependence on $\alpha$ is continuous.
\end{cor}
\begin{proof}
Immediate from Proposition \ref{prop:malik} and Corollary \ref{cor:evaine}.
\end{proof}
Figure \ref{fig:slava} gives the
 value of $\delta$ in terms of $\alpha$ and $\Gal(K'/\mathbb{Q})$.
The numerical value of $\delta$ for $\alpha=0.4014\dotsc$ is given
in the rightmost column.
The table was computed by means of the finite-group package GAP \cite{GAP4}
running on a GNU/Linux box.
The sixteen transitive permutation groups on $\{1,\dotsc,6\}$
have been labelled as in \cite{CHM}.
\begin{figure}
\begin{tabular}{lcl}
$\Gal(K'/\mathbb{Q})$ & $\delta$ & for $\alpha = 0.4014\dotsc$\\ \hline
$C(6) $ & $ \frac{1}{6} 2^{4\alpha}$ & $0.5072\dotsc$\\
$D_6(6) $ & $ \frac{1}{6} 2^{4\alpha}$ & $0.5072\dotsc$\\
$D(6) $ & $ \frac{1}{4} 2^{2\alpha}+ \frac{1}{12} 2^{4\alpha}$ & $0.6897\dotsc$\\
$A_4(6) $ & $\frac{1}{4} 2^{2\alpha} + \frac{1}{12} 2^{4 \alpha}$ & $0.6897\dotsc$\\
$F_{18}(6) $ & $ \frac{2}{9} 2^{2 \alpha} + \frac{1}{18} 2^{4\alpha}$ & $0.5567\dotsc$\\
$2A_4(6) $ & $ \frac{1}{8} 2^{2\alpha}+ \frac{1}{8} 2^{3\alpha}+ \frac{1}{24} 2^{4\alpha}$ & $0.6328\dotsc$\\
$S_4(6d) $ & $\frac{3}{8} 2^{2\alpha}+ \frac{1}{24} 2^{4\alpha}$ & $0.7809\dotsc$\\
$S_4(6c) $ & $\frac{1}{4} 2^{\alpha}+ \frac{1}{8} 2^{2\alpha}+ \frac{1}{24} 2^{4\alpha}$ & $0.6750\dotsc$\\
$F_{18}(6):2 $ & $ \frac{13}{36} 2^{2\alpha}+ \frac{1}{36} 2^{4\alpha}$ & $0.7145\dotsc$\\
$F_{36}(6) $ & $ \frac{13}{36} 2^{2\alpha}+ \frac{1}{36} 2^{4\alpha}$ & $0.7145\dotsc$\\
$2S_4(6) $ & $\frac{1}{8} 2^{\alpha}+ \frac{3}{16} 2^{2\alpha}+ \frac{1}{16} 2^{3\alpha} + \frac{1}{48} 2^{4 \alpha}$ & $0.6996\dotsc$\\
$L(6) $ & $ \frac{2}{5}+ \frac{1}{4} 2^{2\alpha}+  \frac{1}{60} 2^{4\alpha}$ & $0.8868\dotsc$\\
$F_{36}(6):2 $ & $\frac{1}{6} 2^{\alpha}+ \frac{13}{72} 2^{2\alpha}+ \frac{1}{12} 2^{3\alpha}+ \frac{1}{72} 2^{4\alpha}$ & $0.7693\dotsc$\\
$L(6):2 $ & $ \frac{1}{5}+ \frac{1}{4} 2^{\alpha}+ \frac{1}{8} 2^{2\alpha}+ \frac{1}{120} 2^{4\alpha}$ & $0.7736\dotsc$\\
$A_6 $ & $ \frac{2}{5}+ \frac{17}{72} 2^{2\alpha}+ \frac{1}{360} 2^{4\alpha}$ & $0.8203\dotsc$\\
$S_6 $ & $ \frac{1}{5}+ \frac{7}{24} 2^{\alpha}+ \frac{17}{144} 2^{2\alpha}+ \frac{1}{48} 2^{3\alpha}+ \frac{1}{720} 2^{4\alpha}$ & $0.8434\dotsc$\\
\end{tabular}
\caption{Galois groups and corresponding averages}\label{fig:slava}
\end{figure}
\subsection{Rational points and ranks of Jacobians}
\begin{lem}\label{lem:living}
Let $K$ be a number field. Let $f\in K\lbrack x\rbrack$ be an irreducible
polynomial of degree $6$. Let $C$ be the curve $y^2 = f(x)$. For every
$d\in K^*$, let $C_d$ be the curve $d y^2 = f(x)$. Let $c$ be a divisor
of the Jacobian $J_C$ of $C$. Then the image of $C_d(K)$ under the
composition
\begin{equation}\label{eq:phidc} \phi_{d,c}:C_d(K)
\xrightarrow{(x,y)\mapsto (x,y \sqrt{d})} C 
\xrightarrow{P\mapsto \Cl(P) - c} J_C\end{equation}
generates an abelian group of rank at most $\rnk(J_{C_d}(K))+1$.
\end{lem}
\begin{proof}
We have a commutative diagram
\[ \begin{CD} 
       C_d(K) @>{\mu}>> C\\
@VV V @VV V \\
J_{C_d} @>{\sim}>> J_C,\\
       \end{CD} 
       \] 
where the equivalence at the bottom is induced by 
$t_d:(x,y)\mapsto (x,y \sqrt{d})$, and the map on the left is given by
$P\to \Cl(P) - t_d^*c$. If $C_d(K)$ is empty, we have nothing to prove.
Assume that $C_d(K)$ has at least one element, and call it $P_0$. The
image of 
\[C_d(K) \xrightarrow{P\mapsto \Cl(P) - \Cl(P_0)} J_{C_d}\]
is an abelian group of rank no greater than $\rnk(J_{C_d(K)})$. Hence,
the corresponding subgroup of $J_C$ under the isomorphism $J_{C_d}
\equiv J_C$ has rank at most $\rnk(J_{C_d(K)})$ as well. This subgroup
is the image of $C_d(K)$ under the composition
\[C_d(K)\xrightarrow{(x,y)\mapsto (x,y \sqrt{d})} C 
\xrightarrow{P\mapsto \Cl(P) - (t_d^{-1})^*(\Cl(P_0))} J_C .\]
If we displace the subgroup by $(t_d^{-1})^*(\Cl(P_0)) - c$, 
we obtain $\phi_{d,c}(C_d(K))$. The statement follows immediately.
\end{proof}
\begin{prop}[Mumford]\label{prop:mumf} Let $K$ be a number field. Let $C$ be a complete
non-singular curve over $K$ with genus $g\geq 2$. Let $J_C$ be the 
Jacobian of $C$. Let $h_{\delta}$ be the height function on $J\times J$
induced by the theta divisor; let $P\to \Cl(P) - c_0$ be a normalized
embedding of $C$ in $J$, where $c_0$ is an appropiately chosen
divisor of $C$. (See \cite{La}, p.\,113 and p.\,120.)
Define
\begin{equation}\label{eq:braks} \langle P_1,P_2\rangle = 
- h_{\delta}(P_1,P_2),\;\;
|P| = \sqrt{\langle P,P\rangle} .\end{equation}
Then, for any $P_1, P_2 \in C(\overline{K})$, $P_1\ne P_2$,
\[2 g \langle P_1,P_2\rangle \leq |P_1|^2 + |P_2|^2 + O(1) .\]
Notice moreover that
\[\langle P,P\rangle = 2 g h_{c_0}(P) + O(1)\]
for any point $P$ of $C$. The implied constants depend only on $C$ and
on the choice of $c_0$.
\end{prop}
\begin{proof}
See \cite{La}, Thm. 5.10 and Thm. 5.11. (The original formulation in \cite{Mu}
is restricted to rational points.)
\end{proof}
\begin{lem}\label{lem:socmu}
Let $f\in \mathbb{Z}\lbrack x\rbrack$ be an irreducible polynomial of 
degree $6$. For every square-free integer $d$, let $C_d$ be the curve
$d y^2 = f(x)$. Let $c_0$ be as in Proposition \ref{prop:mumf}
and let $\phi_{d,c_0}$ be as in (\ref{eq:phidc}).
 Then there is a constant $\kappa_f$ such that, if $d$ is a square-free
integer $d$ larger than $\kappa_f$, there are at most
$8$ points $P$ in $C_d(\mathbb{Q})$ for which $\phi_{d,c_0}(P)$ is torsion.
\end{lem}
\begin{proof}
Consider the divisor $h_y$ on $C=C_1$. It is clear that
$P = (x, y \sqrt{d})$ lies on $C$, and that $h_y(P) = y \sqrt{d} + O_f(1)$. 
 By \cite{La}, Ch. 4, Cor. 3.5,
\[(1-\epsilon) \kappa_0 h_{c_0}(P) - \kappa_{\epsilon} \leq 
h_y((x,y \sqrt{d})) \leq (1 + \epsilon) \kappa_0 h_{c_0}(P) + \kappa_{
\epsilon}\]
for every $\epsilon > 0$ and some $\kappa_0$, $\kappa_{\epsilon}$, where
$\kappa_{\epsilon}$ depends only on $f$ and $\epsilon$. By 
\cite{La}, Thm. 5.10, $h_{c_0}(P) = \frac{1}{4} |P|^2 + O(1)$,
where $|P|$ is as in (\ref{eq:braks}). Suppose $\phi_{d,c_0}(P)$
is torsion. Then $|P|=0$. We obtain that $h_{c_0}(P)$ is bounded,
and hence so is $h_y(P)$; yet, by Lemma \ref{lem:pseusi}, we know
that either $h_y(P) \geq \frac{3}{8} \log |d| + O_f(1)$, or $y=0$,
or $(x,y \sqrt{d})$ is one of the two points of $C$ at infinity.
If $d$ is large enough, we can conclude that either $y=0$ 
or $(x,y \sqrt{d}) = \pm \infty$. There are six solutions to the
former equation and two to the latter.
\end{proof}
As is usual, we speak of the rank of a curve $C$ when we mean the
rank of its Jacobian. Thus $\rnk(C)$ and $\rnk(J_C)$ are the same.
\begin{prop}\label{prop:ohcan}
Let $f\in \mathbb{Z}\lbrack x\rbrack$ be an irreducible polynomial
of degree $6$. For every square-free integer $d$, let $C_d$ be the
curve $d y^2 = f(x)$. 
For $0< \epsilon < 1/2$, let
\[S = \{(x,y)\in C_d(\mathbb{Q}) : (1-\epsilon) x_0 \leq H(x) 
\leq (1+\epsilon) x_0\}.\]
Then
\[\# S \ll_{\epsilon} 
2^{(\alpha + O_{\epsilon}(x_0^{-2})) \rnk(C_d(K))} ,\]
where $\alpha$ is as in (\ref{eq:calor}) and the implied constants
in $\ll$ and $O_{\epsilon}(x_0^{-2})$ 
depend only on $f$ and $\epsilon$.
\end{prop}
\begin{proof}
Let $C=C_1$.
Let $c_0$ be as in Proposition \ref{prop:mumf},
$\phi_{d,c_0}$ as in (\ref{eq:phidc}). 
By Lemma \ref{lem:living}, all points of $\phi_{d,c_0}(C_d(\mathbb{Q}))$
lie on an abelian subgroup of $J_C$ of rank at most $\rnk (J_{C_d}(K)) + 1$.
The function \[\langle P_1,P_2\rangle = - h_{\delta} : J\times J \to 
\mathbb{R}\]
in Proposition  \ref{prop:mumf} is a positive, symmetric quadratic form. 
We have $\langle P,P\rangle = 0$ if and only if $P$ is torsion. 
Hence there is a map \[\iota : \phi_{d,c_0}(C_d(\mathbb{Q})) \to
\mathbb{R}^{\rnk(J_{C_d}(K)) + 1}\] with torsion kernel such that
 the inner product in $\mathbb{R}^{\rnk(J_{C_d}(K))+1}$ is taken 
 back to the inner product $\langle \cdot,\cdot\rangle$ in
$J$.
By Lemma
\ref{lem:socmu}, there are at most $8$ points $P$ in $C_d(\mathbb{Q})$
such that $\phi_{d,c_0}(P)$ is torsion. Thus the kernel of $\iota$
has cardinality at most $8$, and we can focus on determining the cardinality
of $\iota(\phi_{d,c_0}(S))$.  

Let $P = (x,y) \in S$. 
By \cite{La}, Cor. IV.3.5 and Thm. V.5.10, 
\[(1 - \epsilon)^2 \kappa_0 x_0 - \kappa_\epsilon \leq |\phi_{d,c_0}(P)|^2 \leq
(1 + \epsilon)^2 \kappa_0 x_0 + \kappa_{\epsilon} \]
for some $\kappa_0$, $\kappa_{\epsilon}$, the latter depending on $\epsilon$.
For any $P_1, P_2 \in C_d(\mathbb{Q})$, $P_1\ne P_2$,
 we have by Proposition \ref{prop:mumf} that
\[\langle \phi_{d,c_0}(P_1),
\phi_{d,c_0}(P_2)\rangle \leq \frac{|\phi_{d,c_0}(P_1)|^2 + 
|\phi_{d,c_0}(P_2)|^2}{4} + O(1) .\]
Hence, for $P_1 = (x_1,y_1)$, $P_2 = (x_2,y_2)$ in $S$, $P_1\ne P_2$, 
\[\begin{aligned}
\frac{\langle \iota(\phi_{d,c_0}(P_1)), \iota(\phi_{d,c_0}(P_2))\rangle}{
|\iota(\phi_{d,c_0}(P_1))| |\iota(\phi_{d,c_0}(P_2))|} &=
\frac{\langle \phi_{d,c_0}(P_1),\phi_{d,c_0}(P_2)\rangle}{
|\phi_{d,c_0}(P_1)| |\phi_{d,c_0}(P_2)|} \\ &\leq
\frac{\frac{|\phi_{d,c_0}(P_1)|}{|\phi_{d,c_0}(P_2)|} + 
\frac{|\phi_{d,c_0}(P_2)|}{|\phi_{d,c_0}(P_1)|} 
+ O(x_0^{-2})}{4} \\ &\leq
\frac{\frac{(1 - \epsilon)^2}{(1 + \epsilon)^2} + 
 \frac{(1 + \epsilon)^2}{(1 - \epsilon)^2} + O_{\epsilon}(x_0^{-2})}{4}
\\&= \frac{(1 + 6 \epsilon^2 + \epsilon^4)/(1 - \epsilon^2)^2}{2} +
O_{\epsilon}(x_0^{-2}) .\end{aligned}\]
We are in the same situation as in subsection \ref{subs:kiss}: we
have to bound the number of points that can lie in $\mathbb{R}^n$
with an angle of at least $60^\circ - \epsilon$ between any two of them.
As before, we can apply (\ref{eq:calor}). We obtain 
\[\# \iota(S) \ll 2^{(\alpha + O(\epsilon^2) + O_{\epsilon}(x_0^{-2})) 
(\rnk(J_{C_d}(K)) + 1)} ,\]
from which the statement immediately follows.
\end{proof}
\begin{cor}\label{cor:etern}
 Let $f(x,y) \in \mathbb{Z}\lbrack x,y\rbrack$ be an irreducible
homogeneous polynomial of degree $6$. Let $d$ be a square-free integer. 
Consider the curve $C_d : d y^2 = f(x,1)$. Let $S$ be the set
\[\{(x,y,z)\in \mathbb{Z}^3 : N^{1 - \epsilon} \leq |x|, |z| \leq
N^{1 + \epsilon}, d y^2 = f(x,z)\} .\]
Then
\[\# S \ll_{\epsilon} 2^{(\alpha + O_{\epsilon}((\log N)^{-2})) \rnk(C_d),}\]
where $\alpha$ is as in (\ref{eq:calor}) and the implied constants
in $\ll$ and $O_{\epsilon} (x_0^{-2})$ 
depend only on $f$ and $\epsilon$.
\end{cor}
\begin{proof}
Immediate from Prop. \ref{prop:ohcan}.
\end{proof}
\subsection{The square-free sieve for homogeneous sextics}
\begin{prop}\label{prop:bigsix}
Let $f\in \mathbb{Z}\lbrack x,z\rbrack$ be an irreducible homogeneous 
polynomial of degree $6$. Let
\[\delta(N) = \{-N\leq x,z\leq N : \gcd(x,y)=1, \exists p > N^{1/2} \;\st\; p^2|f(x,z)\} .\]
Then, for every $\epsilon>0$,
\begin{equation}\label{eq:pathe}\delta(N) \ll X (\log X)^{-\beta+\epsilon},
\end{equation}
where $\beta$ is given in Figure \ref{fig:gerh} in terms of the Galois 
group of the splitting field of $f(x,1)=0$. The implied constant
depends only on $f$ and $\epsilon$.
\end{prop}
\begin{proof}
Let $A = \max_{1\leq x,z\leq N} f(x,z)$. Clearly $A\ll N^6$. We can
write
\[\begin{aligned}
\delta(N) &\leq \sum_{N^{1/2}<p<\sqrt{A/M}} 
\# \{-N\leq x,z\leq N : p^2 | f(x)\}\\
&+ \{-N^{1-\epsilon}\leq x,z\leq N^{1-\epsilon} : \gcd(x,z)=1, \exists p>N^{1/2} \text{ s.t. }
p^2 | f(x,z) \} \\
&+ \mathop{\sum_{1\leq |d|\leq M}}_{\text{$d$ square-free}} 
\#\{x,z\in \mathbb{Z}^2 : N^{1-\epsilon}\leq 
\max(|x|,|z|) \leq N,
d y^2 = f(x,z) \} .\end{aligned}\]
Let $M\leq N^2$. By Lemma \ref{lem:facil},
the first term is at most
\[\sum_{N^{1/2} < p < \sqrt{A/M}} 3\ll \frac{3 \sqrt{A/M}}{\log \sqrt{A/M}} \ll \frac{N^3 M^{-1/2}}{\log N} .\]
The second term is clearly no greater than $N^{2-\epsilon}$.
It remains to bound
\[\mathop{\sum_{1\leq |d|\leq M}}_{\text{$d$ square-free}} B(d),\]
where
\[B(d) = \# \{x,z\in \mathbb{Z}^2 : N^{1-\epsilon}\leq 
\max(|x|,|z|) \leq N,
d y^2 = f(x,z) \} .\]
By Corollary \ref{cor:etern},
\[B(d) \ll 2^{(\alpha + O(\epsilon^2) + O((\log N)^{-2})) \rnk(C_d)} .\]
By (\ref{eq:cru}) and Proposition \ref{prop:malik}, 
we can omit the term $O((\log N)^{-2})$. Notice also that we can
treat $d$ negative just as we do $d$ positive, by working with $-f(x,z)$
instead of $f(x,z)$.
Hence it is left to estimate
\[S(M) = \mathop{\sum_{1\leq d\leq M}}_{\text{$d$ square-free}}
 2^{(\alpha + O(\epsilon^2)) \rnk(C_d)} .\]
Corollary \ref{cor:fasc} gives us
$S(M)\ll M (\log M)^{\delta - 1 + O(\epsilon^2)}$,
where $\delta$ is as in Figure \ref{fig:slava}. Set $M = N^2 (\log N)^{
-\frac{2}{3} \delta}$. Then
\[
S(M)\ll N^2 (\log N)^{\frac{1}{3} \delta - 1 + O(\epsilon^2)},\;\;\;
\frac{N^3 M^{-1/2}}{\log N} 
\ll N^2 (\log N)^{\frac{1}{3} \delta - 1 + O(\epsilon^2)} .\]
\end{proof}
\begin{figure}
\begin{tabular}{l l|l l}
Gal & $\beta$ & Gal & $\beta$ \\ \hline 
$C(6) $ & $0.8309\dotsc$ &$D_6(6) $ & $0.8309\dotsc$\\
$D(6) $ & $0.7700\dotsc$ &$A_4(6) $ & $0.7700\dotsc$\\
$F_{18}(6) $ & $0.8144\dotsc$ &$2A_4(6) $ & $0.7890\dotsc$\\
$S_4(6d) $ & $0.7396\dotsc$ &$S_4(6c) $ & $0.7749\dotsc$\\
$F_{18}(6):2 $ & $0.7618\dotsc$ &$F_{36}(6) $ & $0.7618\dotsc$\\
$2S_4(6) $ & $0.7667\dotsc$ &$L(6) $ & $0.7043\dotsc$\\
$F_36(6):2 $ & $0.7435\dotsc$ &$L(6):2 $ & $0.7421\dotsc$\\
$A_6 $ & $0.7265\dotsc$ &$S_6 $ & $0.7188\dotsc$\\
\end{tabular}
\caption{Square-free sieve for sextics: exponents}\label{fig:gerh}
\end{figure}

\section{Square-free values of polynomials}\label{sec:cuanumer}
We can now state our main unconditional results.
\begin{thm}\label{thm:croiss}
Let $f\in \mathbb{Z}\lbrack x\rbrack$ be an irreducible polynomial
of degree $3$. Then the number of positive integers $x\leq N$ for
which $f(x)$ is square-free is given by
\begin{equation}\label{eq:atrya}
N \prod_p \left(1 - \frac{\ell(p^2)}{p^2}\right) + O(N (\log N)^{-\beta}),
\end{equation}
where
\[\begin{aligned}
\beta &= \begin{cases} 0.8061\dotsc
&\text{ if the discriminant of $f$ is a square,}\\ 
0.6829\dotsc
&\text{ if the discriminant of $f$ is not a square,}\end{cases}\\
\ell(m) &= \#\{x\in \mathbb{Z}/m : f(x) \equiv 0 \mo m\}.\end{aligned}\]
The implied constant in (\ref{eq:atrya}) depends only on
$f$.
\end{thm}
\begin{proof}
By Propositions \ref{prop:owerr} and \ref{prop:bethoo}.
\end{proof}
\begin{thm}
Let $f\in \mathbb{Z}\lbrack x,y\rbrack$ be a homogeneous polynomial
of degree no greater than $6$. Then the number of integer pairs
$(x,y)\in \mathbb{Z}^2\cap \lbrack -N,N\rbrack^2$
 for which 
$f(x,y)$ is square-free is given by
\[4 N^2 \prod_p \left(1 - \frac{\ell_2(p^2)}{p^4}\right) +
\begin{cases} O(N^{4/3} (\log N)) &\text{if $\deg_{\irr}(f) \leq 4$,}\\
O(N^{(5+\sqrt{113})/8 + \epsilon}) 
&\text{if $\deg_{\irr}(f) = 5$,}\\
O(N^2 (\log N)^{-\beta + \epsilon}) 
&\text{if $\deg_{\irr}(f) = 6$,}\end{cases}\]
where $\epsilon$ is an arbitrarily small positive number, 
$\beta$ is given in Figure \ref{fig:gerh} in terms of the Galois 
group of the splitting field of $f(x,1)=0$,
$A$ depends only on $f$, the implied constant
depends only on $f$ and $\epsilon$, $\deg_{\irr}$ denotes the
degree of the irreducible factor of $f$ of largest degree, and
\[\ell_2(m) = \#\{(x,y)\in 
(\mathbb{Z}/m)^2 : f(x,y) \equiv 0 \mo m\}.\]
\end{thm}
\begin{proof}
By Propositions \ref{prop:aogi}, \ref{prop:aulait}, \ref{prop:fourcoffees}
\ref{prop:jocha} and \ref{prop:bigsix}. In the case of $\deg f = 1$
and $\deg f = 2$, we reset $M = C N^2$ within Prop. \ref{prop:aogi},
where $C$ is any constant such that $|f(x,y)|\leq C N^2$ for all
$1\leq x,y\leq N$.
\end{proof}
In the same way as above, we can obtain unconditional results from
the propositions in subsection \ref{subs:avmul} by applying Propositions
\ref{prop:bethoo}, \ref{prop:aulait}, \ref{prop:fourcoffees},
\ref{prop:jocha} and \ref{prop:bigsix}.
\section{Previous work and work to do}\label{sec:bestpast}
The best bounds known before now are listed in the first table of the
introduction. The case $\deg_{\irr}(P(x))=2$ was dealt with by
Estermann \cite{Es}, the case $\deg_{\irr}(P(x))=3$ by Hooley 
(\cite{Ho}, Ch. IV). All entries for $P(x,y)$ homogeneous,
$3\leq \deg_{\irr}(P(x))\leq 6$, are due to Greaves \cite{Gr}.
The bound given in the main theorem of \cite{Gr} for $\deg_{\irr} = 6$
is actually $N^2/(\log N)^{1/3}$, rather than $N^2/(\log N)^{1/2}$.
Ramsay (\cite{Ra}, 1991) showed that a slight amendment is sufficient
to improve the exponent from $(\log N)^{1/3}$ to $(\log N)^{1/2}$.
In our formulation, it is enough to substitute Lemma \ref{lem:ramsay}
for \cite{Gr}, Lemma 1, within the proof of \cite{Gr}, Lemma 2.

Most of the results in the present paper should carry over fairly readily
to polynomials over number fields. The estimates in \S \ref{sec:gust},
coming from Diophantine geometry, can be generalized at least as easily
as bounds coming from sieves. However, we have stated the final results
only over $\mathbb{Q}$, except for the cases in which a general treatment
required little or no additional space. The main reason for this limitation
is that we do not know a priori what kind of generalization is desirable.
Given a polynomial $P\in \mathscr{O}_K\lbrack x\rbrack$, should we examine
only its values for $x\in \mathbb{Z}$, or should we let $x$ range over all
of $\mathscr{O}_K$? In the latter case, how do we order $\mathscr{O}_K$?
A simple ordering by norm will not do, as there are infinitely many
elements of $\mathscr{O}_K$ of norm below any given $c>1$, unless
$K$ is imaginary quadratic. One may choose a basis of $\mathscr{O}_K$
and use this basis to inject a box $\{1,2,\dotsc,N\}^{\deg K}$ into
$\mathscr{O}_K$, but this procedure is neither canonical nor necessarily
natural. It is probably best to wait to see what will be demanded by
applications, and to hope that abstract statements such as Prop. \ref{prop:ridd}
and Cor. \ref{cor:yugo} 
will accommodate the required change in the objects of study.

\section{Acknowledgments}
I am indebted to my advisor, Henryk Iwaniec, for his guidance. Thanks are
due as well to Jordan Ellenberg and Joseph Silverman, for their advice,
and to Keith Ramsay, for our discussions on his unpublished work.

\appendix
\section{Lattices and solutions}\label{sec:appa}
We write $\Disc P$ for the discriminant of a polynomial $P\in \mathscr{O}_K\lbrack x\rbrack$.
If $P$ is square-free, then 
$\Disc(P)$ is a non-zero element of $\mathscr{O}_K$, and
\[\gcd(P(x),P'(x))|\Disc(P)\] for every $x\in \mathscr{O}_K$.

\begin{lem}\label{lem:sols}
Let $K$ be a $\mathfrak{p}$-adic field.
Let $P\in \mathscr{O}_K\lbrack x \rbrack$ be a square-free
polynomial. Then
\[P(x)\equiv 0 \mo \mathfrak{p}^n\]
has at most $\max(|\Disc P|_\mathfrak{p}^{-1} \cdot \deg P, 
|\Disc P|_\mathfrak{p}^{-3})$ roots in $\mathscr{O}_K/\mathfrak{p}^n$.
\end{lem}
\begin{proof}
Let $\pi$ be a prime element of $K$.
If $P$ is of the form $P = \pi Q$ for some $Q\in \mathscr{O}_K \lbrack x \rbrack$, the statement follows
from the statement for $Q$. Hence we can assume $P$ is not of the form $P = \pi G$. Write
$P = P_1 \cdot P_2 \cdot \dotsb \cdot P_l$, $P_i\in \mathscr{O}_K$, $P_i$ irreducible.

If $n\leq 3 v_\mathfrak{p}(\Disc P)$, there are trivially at most 
$\# (\mathscr{O}_K/\mathfrak{p}^n) = |\mathfrak{p}^n|_\mathfrak{p}^{-1} 
\leq |\Disc P|_\mathfrak{p}^{-3}$ roots.
Assume $n> 3 v_{\mathfrak{p}}(\Disc P)$.
 Let $x$ be a root of 
$P(x)\equiv 0 \mo \mathfrak{p}^n$. Let $P_i$ be a factor for which $v_\mathfrak{p}(P_i(x))$ is maximal. By
\[\begin{aligned}
v_\mathfrak{p}(P'(x)) &= v_\mathfrak{p}(\sum_j P_j'(x) \cdot P_1(x) \cdot \dotsb \widehat{P_j(x)}\dotsb \cdot P_n(x))
\\ &\geq \min_j(v_\mathfrak{p}(P(x)) - v_\mathfrak{p}(P_j(x))),\end{aligned}\]
\[\min(v_\mathfrak{p}(P'(x)),v_\mathfrak{p}(P(x)))\leq 
v_\mathfrak{p}(\Disc P),\]
\[v_\mathfrak{p}(P(x))>v_\mathfrak{p}(\Disc P),\] we have that
\[\min_j(v_\mathfrak{p}(P(x)) - v_\mathfrak{p}(P_j(x))) \leq v_\mathfrak{p}(\Disc P)\] and hence
\[v_\mathfrak{p}(P_i(x)) \geq v_\mathfrak{p}(P(x)) - 
v_\mathfrak{p}(\Disc P)\geq n - v_\mathfrak{p}(\Disc P)\geq 2 v_\mathfrak{p}(\Disc P) + 1.\]
On the other hand $\gcd(P_i(x),P_i'(x))|\Disc P$, and thus 
\[v_\mathfrak{p}(P_i'(x))\leq v_\mathfrak{p}(\Disc P).\]
By Hensel's lemma we can conclude that $P_i$ is linear. Since $v_\mathfrak{p}(P_i(x))\geq n - v_\mathfrak{p}(\Disc P)$,
$x$ is a root of \[P_i(x)\equiv 0 \mo \mathfrak{p}^{n-v_\mathfrak{p}(\Disc P)} .\]
Since $P_i$ is linear and not divisible by $\mathfrak{p}$, it has at most one root in the ring
$\mathscr{O}_K/\mathfrak{p}^{n-v_\mathfrak{p}(\Disc P)}$. There are at most $v_\mathfrak{p}(\Disc P)$ elements of 
$\mathscr{O}_K/\mathfrak{p}^n$ reducing to this root. Summing over all $i$ we obtain that there are at
most $l\cdot v_\mathfrak{p}(\Disc P)$ roots of $P(x)\equiv 0 \mo \mathfrak{p}^n$ in $\mathbb{Z}/\mathfrak{p}^n$. Since
$l\leq \deg P$, the statement follows.
\end{proof}

For every non-zero $\mathfrak{m}\in I_K$, we define $h(\mathfrak{m})$ to
be the positive integer generating $\mathfrak{m} \cap \mathbb{Z}$.

\begin{lem}\label{lem:dontaskmemyadvice}
 Let $K$ be a number field. Let $\mathfrak{m}$ be a non-zero ideal
of $\mathscr{O}_K$. Let $P\in \mathscr{O}_K\lbrack x\rbrack$ be a square-free
polynomial. Then 
\[\{x\in \mathbb{Z} : P(x)\equiv 0 \mo \mathfrak{m}\}\]
is the union of at most $|\Disc P|^3 \tau_{\deg P}(\rad(h(\mathfrak{m}))))$
arithmetic progressions of modulus $h(\mathfrak{m})$. 
\end{lem}
\begin{proof}
  By Lemma \ref{lem:sols}, for every $\mathfrak{p}|\mathfrak{m}$,
the equation
\[P(x)\equiv 0 \mo \mathfrak{p}^n\]
has at most 
$|\Disc P|_\mathfrak{p}^{-3} \deg P$ roots in $\mathscr{O}_K/\mathfrak{p}^n$.
For any ideal $\mathfrak{a}$,
the intersection of $\mathbb{Z}$ with a set of the form
\[\{x\in \mathscr{O}_K : x\equiv x_0 \mo \mathfrak{a}\}\]
is either the empty set or an arithmetic progression of modulus 
$h(\mathfrak{a})$. This is in particular true for $\mathfrak{a} = \mathfrak{p}^n$; the set 
\[\{x\in \mathbb{Z} : x\equiv x_0 \mo \mathfrak{p}^n\}\]
is the union of at most 
$|\Disc P|_\mathfrak{p}^{-3} \deg P$ arithmetic progressions of modulus
$h(\mathfrak{p}^n)$. 

Now consider a rational prime $p$ at least one of whose prime ideal
divisors divides $m$. Write $m = \mathfrak{p}_1^{n_1} \mathfrak{p}_2^{n_2}
\dotsb \mathfrak{p}_k^{n_k} \mathfrak{m}_0$, where $\mathfrak{p}_1,\dotsc,
\mathfrak{p}_k | p$, $\gcd(\mathfrak{m}_0,p)=1$, 
$\frac{n_1}{e_1}\geq \frac{n_2}{e_2}\geq \dotsb \geq \frac{n_k}{e_k}$ and
$e_1,\dotsc, e_k$ are the ramification degrees of $\mathfrak{p}_1,\dotsc,
\mathfrak{p}_k$. The set
\[\{x\in \mathbb{Z} : x\equiv x_0 \mo \mathfrak{p}_1^{n_1} \dotsb
\mathfrak{p}_k^{n_k}\}\]
is the intersection of the sets
\[\{x\in \mathbb{Z} : x\equiv x_0 \mo \mathfrak{p}_j^{n_j}\},
\text{\;\; $1\leq j\leq k$.}
\]
At the same time, it is a disjoint union of arithmetic progressions
of modulus \[h(\mathfrak{p}_1^{n_1}\dotsb\mathfrak{p}_k^{n_k}) =
h(\mathfrak{p}_1^{n_1}).\] Since
\[\{x\in \mathbb{Z} : x\equiv x_0 \mo \mathfrak{p}_1^{n_1}\}
\] is the disjoint union of at most
$|\Disc P|_\mathfrak{p}^{-3} \deg P$ arithmetic progressions
of modulus $h(\mathfrak{p}_1^{n_1})$, it follows that
\[\{x\in \mathbb{Z} : x\equiv x_0 \mo \mathfrak{p}_1^{n_1} \dotsb
\mathfrak{p}_k^{n_k}\}\]
is the disjoint union of at most
$|\Disc P|_\mathfrak{p}^{-3} \deg P$ arithmetic progressions
of modulus $h(\mathfrak{p}_1^{n_1}) = h(\mathfrak{p}_1^{n_1} \dotsc \mathfrak{p}_k^{n_k})$. The statement follows.
\end{proof}

\begin{lem}\label{lem:sollat} 
Let $K$ be a number field. 
Let $\mathfrak{m}$ be a non-zero ideal of $\mathscr{O}_K$.
Let $P\in \mathscr{O}_K\lbrack x,y\rbrack $ be a non-constant and square-free
 homogeneous
polynomial. Then the set
\[S = \{(x,y)\in \mathbb{Z}^2 : \gcd(x,y) = 1, \mathfrak{m}|P(x,y)\}\]
is the union of at most 
$|\Disc P|^3\cdot \tau_{2\deg P}(\rad(h(\mathfrak{m})))$ 
disjoint sets of the form
$L \cap \{(x,y)\in \mathbb{Z}^2 : \gcd(x,y) = 1\},$
$L$ a lattice of index $\lbrack \mathbb{Z}^2 : L\rbrack = h(\mathfrak{m})$.
\end{lem}
\begin{proof}
Let $\mathfrak{p}|\mathfrak{m}$. Let $n = v_\mathfrak{p}(\mathfrak{m})$.
Let $r_1,r_2,\dotsb,r_k\in \mathscr{O}_K/\mathfrak{p}^n$ be the roots of
$P(r,1)\cong 0 \mo \mathfrak{p}^n$. Let $r'_1,r'_2,\dotsb,r'_{k'}\in
\mathscr{O}_K/\mathfrak{p}^n$ be such roots of
$P(1,r)\cong 0 \mo \mathfrak{p}^n$ as satisfy $\mathfrak{p}|r$.
 Then the set of solutions to
$P(x,y)\cong 0 \mo \mathfrak{p}^n$ in 
\[\{(x,y)\in \mathbb{Z}^2: \mathfrak{p}\nmid \gcd(x,y)\}\] is the union of the disjoint sets
\[\{(x,y)\in \mathbb{Z}^2: \mathfrak{p}\nmid \gcd(x,y), x\equiv r_i y \mo \mathfrak{p}^n\},\]
\[\{(x,y)\in \mathbb{Z}^2: \mathfrak{p}\nmid \gcd(x,y), y\equiv r_i x \mo \mathfrak{p}^n\} .\] 
The rest of the argument is as in Lemma 
\ref{lem:dontaskmemyadvice}.\end{proof}

\begin{lem}\label{lem:penult}
Let $S\subset \mathbb{R}^2$ be a sector.
Let $L$ be a lattice not contained in 
$\{(x,y)\in \mathbb{Z}^2: \gcd(x,y)\ne 1\}$.
Then
\[\# (\{-N\leq x,y\leq N: \gcd(x,y)=1\} \cap S \cap L)\]
equals \[\frac{\Area(S \cap \lbrack - N , N \rbrack^2)}{
                    \lbrack \mathbb{Z}^2 : L\rbrack}
 \prod_p \left(1-\frac{1}{p^2} \right) + O(N \log N) .\]
The implied constant is absolute.
\end{lem}
\begin{proof}
We can write 
$\#(\{1\leq x,y\leq N: \gcd(x,y)=1\} \cap S \cap L)$ as
\[ \sum_{m=1}^N \mu(m) \cdot
\# (\{1\leq x,y\leq N: m|x, m|y\} \cap S \cap L) .\]
Since $L$ is not of the form $\ell L'$, $\ell>1$, $L'\subset \mathbb{Z}^2$
 a lattice, we have
that $\frac{1}{m} (L \cap m \mathbb{Z}^2)$ is a lattice of index
$\lbrack \mathbb{Z}^2 : \frac{1}{m} (L\cap m \mathbb{Z}^2)\rbrack = 
\lbrack \mathbb{Z}^2 : L\rbrack$. Thus
\[\# (\{1\leq x,y\leq N: m|x, m|y\} \cap S \cap L) = 
A \cdot 
N^2/m^2 \lbrack \mathbb{Z}^2 : L\rbrack + O(N/m) ,\]
where $A = \Area(S \cap \lbrack -1, 1 \rbrack^2)$. Hence
\[\begin{aligned}
\# &(\{1\leq x,y,\leq N: \gcd(x,y)=1\} \cap S \cap L) \\ &= A \cdot
\sum_{m=1}^N \mu(m) N^2/m^2 \lbrack \mathbb{Z}^2 : L\rbrack + O(N \log N),\\
&= \frac{\Area(S \cap \lbrack -N ,N \rbrack^2)}{
\lbrack \mathbb{Z}^2 : L\rbrack} \prod_p \left(1 - \frac{1}{p^2}
\right) + O(N \log N) .
\end{aligned}\]
\end{proof}
The following lemma is better than trivial estimates when 
$L$ is a lattice of index greater than $N$.
\begin{lem}\label{lem:ramsay}
Let $L$ be a lattice. Then
\[\#(\{-N\leq x,y\leq N : \gcd(x,y)=1\} \cap L) \ll 
\frac{N^2}{\lbrack \mathbb{Z}^2 : L\rbrack} + 1.\]
\end{lem}
\begin{proof}
Let $M_0 = \min_{(x,y)\in L} \max(|x|,|y|)$.
By \cite{Gr}, Lemma 1, 
\[\#(\lbrack -N,N\rbrack^2 \cap L) \leq \frac{4 N^2}{\lbrack \mathbb{Z}^2 :
L\rbrack} + O\left(\frac{N}{M_0}\right).\]
If $M_0\geq 
\frac{\lbrack \mathbb{Z}^2 : L\rbrack}{2 N}$ we are done. Assume
$M_0< \frac{\lbrack \mathbb{Z}^2 : L \rbrack}{2 N}$. 
Suppose 
\[\#(\{-N\leq x,y\leq N : \gcd(x,y)=1\} \cap L) > 2.\]
Let $(x_0,y_0)$ be a point
such that $\max(|x_0|,|y_0|) = M_0$. Let $(x_1,y_1)$ be a point
in $\#(\{-N\leq x,y\leq N : \gcd(x,y)=1\} \cap L)$ other than
$(x_0,y_0)$ and $(-x_0,-y_0)$. Since $\gcd(x_0,y_0)=\gcd(x_1,y_1)=1$,
it cannot happen that $0$, $(x_0,y_0)$ and $(x_1,y_1)$ lie on the same
line. Therefore we have a non-degenerate parallelogram
$(0,(x_0,y_0),(x_1,y_1),(x_0+x_1,y_0+y_1))$ whose area has to be at least
$\lbrack \mathbb{Z}^2 : L\rbrack$. On the other hand, its area can
be at most $\sqrt{x_0^2+y_0^2} \cdot \sqrt{x_1^2+y_1^2} \leq 
\sqrt{2} M_0 \cdot \sqrt{2} N = 2 M_0 N$. Since we have assumed
$M_0< \frac{\lbrack \mathbb{Z}^2 : L \rbrack}{2 N}$ we arrive at a
contradiction. \end{proof}

\section{Divisor sums}
Throughout the present paper, we use repeatedly, and sometimes without
mention, the most common
bounds on divisor functions: $\tau(n) \ll n^{\epsilon}$,
$\sum_{n=1}^N \tau_k(n) \ll N (\log N)^{k-1}$,
$\sum_{n=1}^N \omega(n) \ll N \log \log N$, and so on. We also
need the following auxiliary results, which are elementary but not
quite standard.
\begin{lem}\label{lem:crudel}
For every $\epsilon>0$,
\[\sum_{\text{$n$ square-free}} 
\prod_{p|n} \frac{\log X}{p \log p} \ll
X^{\epsilon} ,\]
where the implied constant depends only on $\epsilon$.
\end{lem}
\begin{proof}
For every positive $k$,
\[\begin{aligned}
\sum_{p>k} \frac{1}{p \log p} &\leq \frac{1}{\log k}
  \left(\sum_{k\leq p < k^2} \frac{1}{p} + \frac{1}{2}
\sum_{k^2\leq p < k^4} \frac{1}{p}
 + \frac{1}{4}
\sum_{k^4\leq p < k^8} \frac{1}{p} + \dotsb \right) \\
 &= \frac{\log \log k^2 - \log \log k + O(1/(\log k))}{\log k} \\
 &+ \frac{\log \log k^4 - \log \log k^2 + O(1/(\log k))}{2 \log k} + \dotsb\\
 &= \frac{2 (\log 2 + O(1/(\log k)))}{\log k} .\end{aligned}\]
Hence
\[\prod_{p>k} \left(1 + \frac{1}{p \log p}\right) \ll e^{2 \log 2/\log k
+ O(1/(\log k)^2)} .\]
Now
\[\begin{aligned}
\sum_{\text{$n$ square-free}} 
\prod_{p|n} \frac{\log X}{p \log p} 
&= \prod_{p\leq k} \left(1 + \frac{\log X}{p \log p}\right) 
\prod_{p>k} \left(1 + \frac{\log X}{p \log p}\right) \\
&\ll \prod_{p\leq k} \log X
\prod_{p>k} \left(1 + \frac{1}{p \log p}\right)^{\log X}\end{aligned}\]
for $N>2$. Therefore
\[\begin{aligned}
\mathop{\sum_{1\leq n\leq N}}_{\text{$n$ square-free}} 
\prod_{p|n} \frac{\log n}{p \log p} &\ll
e^{2 \log 2 \log X/\log k
+ O(\log X/(\log k)^2)}
\prod_{p\leq k} \log X \\ &\ll
e^{\frac{2 \log 2 \log N}{\log k}
+ O\left(\frac{\log X}{(\log k)^2}\right)}
e^{\frac{k \log \log X}{\log k} +
O\left(\frac{k \log \log X}{(\log k)^2}\right)}.\end{aligned}\]
Set $k = \log X/ \log \log X$. We obtain
\[\mathop{\sum_{1\leq n\leq N}}_{\text{$n$ square-free}}
\prod_{p|n} \frac{\log X}{p \log p} \ll
 e^{(2 \log 2 + \epsilon') \log X/\log \log X}\]
for every $\epsilon'>0$, where the implied constant depends only
on $\epsilon'$.
\end{proof} 
\begin{lem}\label{lem:grogor}
Let $c>0$ be given. For every $\epsilon>0$,
\[\frac{1}{N} \mathop{\sum_{1\leq n\leq N}}_{\text{$n$ square-free}}
 \tau(n) \prod_{p|n} \frac{c_0 \log N}{\log p} \ll N^{\epsilon},\]
where the implied constant depends only on $c_0$ and $\epsilon$.
\end{lem}
\begin{proof}
Clearly
\[\frac{1}{N} \mathop{\sum_{1\leq n\leq N}}_{\text{$n$ square-free}}
 \tau(n) \prod_{p|n} \frac{c \log N}{\log p} \leq
\frac{1}{N} \mathop{\sum_{1\leq n\leq N}}_{\text{$n$ square-free}}
 \prod_{p|d} \frac{2 c \log N}{\log p} \ll
\sum_{n=1}^{\infty} \prod_{p|n} \frac{\log N^{2 c}}{p \log p} .\]
Proceed as in Lemma \ref{lem:crudel}. 
\end{proof}

\end{document}